\documentclass{amsart}
\usepackage{amsthm,amssymb,amsmath}

\newtheorem{theorem}{Theorem}[section]
\newtheorem{lemma}[theorem]{Lemma}

\theoremstyle{definition}

\newtheorem{example}[theorem]{Example}

\theoremstyle{remark}
\newtheorem{remark}[theorem]{Remark}

\numberwithin{equation}{section}

\newcommand{\const}{{\rm const}}
\newcommand{\supp}{{\rm supp}}

\newcommand{\esssup}{{\rm esssup}}

\newcommand{\esslim}{{\rm esslim}}

\newcommand{\cl}{\mbox{cl}}

\newcommand{\loc}{{\rm loc}}

\begin{document}

\title
[Linear differential equation with distributions]
{On well-posedness of the linear Cauchy problem with the distributional right-hand side and \\ discontinuous coefficients}

\author[D.~Kinzebulatov]{D. Kinzebulatov}

\address{Department of Mathematics, University of Toronto, Toronto, Ontario, Canada M5S 2E4}

\email{dkinz@math.toronto.edu}



\begin{abstract}
The classical result on well-posedness of Cauchy problem for the linear ordinary differential system with the distributional right-hand side and smooth matrix of coefficients plays fundamental role in many applications of distribution theory to ordinary and partial differential equations.
In the present paper we generalize this result to the case of system 
\begin{equation}
\label{abseq1}
x^{\prime}-A(t)x=f,
\end{equation}
where $f$ is a distribution, $A \in \mathbb G^\infty$, $\mathbb G^\infty$ is the space of functions possessing at most first-kind discontinuities together with all their derivatives defined almost everywhere. The left-hand side of system (\ref{abseq1}) contains the product of a distribution and, in general, a discontinuous function, which is undefined in the classical space $\mathcal D'$ of distributions with the smooth test functions. As a result, the Cauchy problem for (\ref{abseq1}) in general has no solution in $\mathcal D'$. In what follows, we consider system (\ref{abseq1}) in the space $\mathcal R'$ of distributions with $\mathbb G^\infty$-test functions, whose elements admit continuous multiplication by functions in $\mathbb G^\infty$, and show that there exists the unique solution of the Cauchy problem for (\ref{abseq1}) which depends continuously on $f$. 
The proof of this result requires investigation of structure of the kernel of operator of restriction of distributions from $\mathcal R$ to $\mathcal D$, of properties of operation of multiplication and of properties of multi-valued (yet, in a sense, continuous) operation of differentiation in $\mathcal R'$.
\end{abstract}

\maketitle

\section{Introduction}
The theory of distributions, which made possible to tackle in a systematic and mathematically rigorous way the discontinuous solutions of ordinary
and partial differential equations, appeared in its finished form in 1950 in the famous monograph by L.~Schwartz \cite{Scw} that contained, among other fundamental results of the new theory, the result on well-posedness of the Cauchy problem
for the linear differential system with the smooth matrix of coefficients and distributional right-hand side. It is the generalization of this result to the case of system
\begin{equation}
\label{ieq1}
x'-A(t)x=f,
\end{equation}
where $f$ is a distribution and $A$ may have at most first-kind discontinuities together with all its derivatives defined almost everywhere, which occupies a central place in the present work.

The linear differential equation (\ref{ieq1}) is an example of an ordinary differential equation whose left-hand side may contain the product of a distribution and a discontinuous function. The study of certain classes of such differential equations started shortly after creation of the distribution theory, and was inspired both by efforts to
extend the field of applicability of the theory of distributions (e.g., the study of linear ordinary differential equations with the distributional coefficients of the order of singularity $\leqslant 1$ in \cite{Kurz1,Kurz2}, also, see \cite{Fil,Tvr1,Tvr2}) as well as by numerous problems of optimal control with the unimodal phase restrictions on control, where the product of a distributional optimal control and the corresponding discontinuous solution of the differential equation arise (see \cite{Mil,Ses,Kinz}). 

The major obstacle on the way of construction of systematic theory of such differential equations in the classical framework of the theory of distributions (by which we mean the theory of the space of linear continuous functionals defined on a certain space of the smooth test functions, e.g., the theory of the space $\mathcal D'$ \cite{Gelf}) -- that is, the impossibility of continuous multiplication of distributions by discontinuous functions (see \cite{DK2,DerKin3}) -- is closely related to the requirement of the smoothness of the test functions (which is crucial for the classical "$(f',\varphi)=-(f,\varphi')$" definition of the distributional derivative \cite{Gelf}). For example, in (\ref{ieq1}) the product $Ax$, determined by the equality (let us consider the scalar case for simplicity)
\begin{equation}
\label{iprod}
(Ax,\varphi)=(x,A\varphi),
\end{equation}
where $\varphi$ is a smooth test function, is undefined, since $A\varphi$ is in general no longer smooth. Formally, there is no need in use of the equality (\ref{iprod}) to define the value of the product $Ax$: there are various definitions of the product $Ax$ in the space $\mathcal D'$ such that $Ax$ coincides with the classical product if $A$ is smooth (see \cite{Sar,Sar2,Sar3}, where the family of distributional products which are invariant under unimodular transformations and satisfy the Leibniz product rule is successfully applied to study of certain classes of ordinary and partial differential equations; see \cite{Tvr2,Fil} and further references therein for other definitions of the product of a distribution and a discontinuous function in the space $\mathcal D'$). Nevertheless, the inevitable lack of continuity in the topology of $\mathcal D'$ of any definition of the product of, in particular, Heaviside function and Dirac delta-function
makes every such definition unacceptable for our proof of the well-posedness of the Cauchy problem for (\ref{ieq1}).
It is the space of distributions with the discontinuous test functions where (\ref{iprod}) defines the continuous operation of multiplication of distributions by discontinuous functions
which can provide proper meaning for the product $Ax$ and, thus, for the system (\ref{ieq1}).

The space of distributions with the test functions that are infinitely differentiable outside the origin where they possess both one-sided limits (together with all their derivatives defined outside the origin) was introduced in \cite{Kur} in application to the problems of construction of self-adjoint operators corresponding to finite-rank perturbations of the $n$-th derivative operator with the support at the origin (also, see \cite{Kur2}). The generalization to the case of test functions of several variables possessing, in certain extent, an arbitrary set of points of discontinuity, was obtained in \cite{DK} in order to provide existence of Nash equilibrium for a class of zero-sum games with first-kind discontinuous payoff functions.
Before we define our space of distributions with discontinuous test functions (in the next section),
let us provide certain heuristics on multiplication of Dirac delta-function by Heaviside function.
Let $\rho \in \mathbb L(\mathbf R)$ be such that $$\supp(\rho) \subset (-1,1) \quad \text{ and  }\quad \int_{\mathbf R} \rho(t)dt=1.$$ We define delta-family $\{\rho_\tau^\varepsilon\}_{\varepsilon>0}$ by the formula
\begin{equation*}
\rho_\tau^\varepsilon(t)=\frac{1}{\varepsilon} \rho\left(\frac{t-\tau}{\varepsilon}\right), \quad t \in \mathbf R,
\end{equation*}
so that $\rho_\tau^\varepsilon \to \delta_\tau$ in $\mathcal D'$ as $\varepsilon \to 0+$. Let $\theta_\tau$ be Heaviside function discontinuous at $\tau$. Then
\begin{equation*}
\theta_\tau \rho_\tau^\varepsilon \to \alpha \delta_\tau
\end{equation*}
in $\mathcal D'$, where the complex coefficient $\alpha \in \mathbf C$, as straightforward calculations show, is given by
\begin{equation}
\label{alpharho}
\alpha=\int_0^\infty \rho(t)dt.
\end{equation}
Thus, in order to avoid multi-valuedness of the product of delta-function and Heaviside function
it is necessary to specify "additional information" on $\delta_\tau$,
which is impossible in the classical space $\mathcal D'$ of distributions with the smooth test functions,
but which is, however, possible in the
space $\mathcal R'$ of distributions defined on the space of discontinuous test functions -- the space of functions which have compact support and which are regulated together with their derivatives of all orders defined almost everywhere (function is called regulated if it possesses both one-sided limits at every point of the interval \cite{Dieu}): there is delta-function defined by the formula
\begin{equation*}
(\delta_\tau^\alpha,\varphi)=\alpha\varphi(\tau+)+(1-\alpha)\varphi(\tau-),
\end{equation*}
so $$\theta_\tau\delta_\tau^\alpha=\alpha\delta_\tau^+,$$ where $(\delta_\tau^+,\varphi)=\varphi(\tau+)$, and $\rho^\varepsilon_\tau \to \delta_\tau^\alpha$ if and only if (\ref{alpharho}) holds. In a similar way we may define the derivatives of the delta-function 
\begin{equation*}
(\delta_\tau^{(k)\alpha},\varphi)=\alpha\varphi^{(k)}(\tau+)+(1-\alpha)\varphi^{(k)}(\tau-),
\end{equation*}
so that $(\rho_\tau^{\varepsilon})^{(k)} \to \delta_\tau^{(k)\alpha}$ if and only if the equality (\ref{alpharho}) holds.

Now, being extending major constructions of the classical distribution space $\mathcal D'$ to the space $\mathcal R'$ of distributions with the discontinuous test functions, one immediately encounters the impossibility to employ the classical definition of the derivative
\begin{equation}
\label{opdif}
(f',\varphi)=-(f,\varphi'),
\end{equation}
even when $\varphi'$ is assumed to be defined almost everywhere: first, in this assumption, (\ref{opdif}) gives rise to a non-linear operator of differentiation, second, as it follows from the above considerations, each distribution possesses whole family of its derivatives, e.g., fix arbitrary $\alpha \in \mathbf C$ and consider delta-family $\{\rho_\tau^\varepsilon\}_{\varepsilon>0}$ which corresponds to a locally absolutely continuous function $\rho$ such that (\ref{alpharho}) holds, to obtain that the delta-function $\delta_\tau^\alpha$ is the derivative of $\theta_\tau$; moreover, if $f$ is a distribution, and $f'$ is its derivative, then $f'+c(\delta_\tau^1-\delta_\tau^0)$ is another derivative of $f$ for every $c$ and every $\tau$ (formally, in \cite{Kur,Kur2} the definition (\ref{opdif}) with $\varphi'$ defined outside the origin is used; however, since the derivatives of the distributions which are considered in \cite{Kur,Kur2} are specified in accordance with the natural multi-valuedness of the differentiation operator, the insufficiencies of the definition (\ref{opdif}) are not crucial for the particular results in \cite{Kur,Kur2}; in \cite{DK} the problem of definition of derivative was in fact avoided, since all elements of the distribution space in \cite{DK} are, in a sense,  measure-type distributions). 

It is natural to require from the definition of the derivative of a distribution to possess the following property of continuity: if $f$ is the distribution, and $g$ is its derivative, then there exists a family $\{f_\varepsilon\}_{\varepsilon>0}$ of locally absolutely continuous functions such that
$$f_\varepsilon \to f, \quad f'_\varepsilon \to g$$ as $\varepsilon \to 0+$, and conversely, if $f_\varepsilon \to f$, and there exists a distribution $g$ such that $f'_\varepsilon \to g$, then the distribution $g$ is the derivative of $f$ (in fact, since this is the only requirement which we impose on the definition of the derivative, this property is \textit{already} a definition of the derivative). \\ [-2mm]

Thus, the aims of the present paper are the following. 

1) To give the (analytical) definition of the derivative which possesses the aforementioned property of continuity.

2) To prove the well-posedness of the Cauchy problem for the linear system 
\begin{equation*}
x'-A(t)x=f,
\end{equation*}
where $f$ is a distribution, and $A$ is regulated together with all its derivatives defined almost everywhere (and, thus, in general discontinuous). \\ [-2mm]

We would like to emphasize the fact that the operation of differentiation in $\mathcal R'$ is multivalued. The formulations of the statements which imply required well-posedness of Cauchy problem (Theorems 3.1 -- 3.5 below) coincide with the formulations of the analogous statements for the space $\mathcal D'$ in \cite{Shi} (see \cite{Scw}), but the aforementioned multi-valuedness of the operation of differentiation requires new (yet natural) definition of the solution as well as totally new proofs.

In fact, the formulation of Theorems 3.1 -- 3.5 allows us to write down the solution of the Cauchy problem explicitly, as the following examples show (the proofs of Examples \ref{example1} and \ref{example2} are provided in Section 3).

\begin{example}
\label{example1}
Let us consider in the space $\mathcal R'$ the following Cauchy problems:
\begin{equation}
\label{pr1}
x'=a\theta_\tau x+b\delta_\tau^\alpha, \quad x=0 \text{ if } t<t_0,
\end{equation}
\begin{equation}
\label{pr2}
x'=a\theta_\tau x+b\delta_\tau^{\prime\alpha}, \quad x=0 \text{ if } t<t_0,
\end{equation}
\begin{equation}
\label{pr3}
x'=a\theta_\tau x+b\delta_\tau^{\prime\prime\alpha}, \quad x=0 \text{ if } t<t_0,
\end{equation}
where $a$, $b \in \mathbf C$, $\tau$, $t_0 \in I$, $t_0<\tau$. Let us define $\zeta_\tau(t)=t-\tau$ for $t>\tau$, $\zeta_\tau(t)=0$ for $t<\tau$.
We call the solution of the differential equation (\ref{pr1}), (\ref{pr2}), (\ref{pr3}) the distribution $x \in \mathcal R'$ which possesses the derivative $x' \in \mathcal R'$ such that after the substitution of $x'$ and $x$ in (\ref{pr1}), respectively, in (\ref{pr2}), (\ref{pr3}), the equation becomes the identity $\mathcal R'$.
Using Theorem 3.2 below, we find that the solution of the Cauchy problem (\ref{pr1}) is the ordinary function given by the formula
\begin{equation*}
x=be^{a\zeta_\tau}\theta_\tau
\end{equation*}
(note that right-hand side of (\ref{pr1}) does not contain the product of a singular distribution and a discontinuous function, and the solution of the Cauchy problem does not depend on $\alpha \in \mathbf C$).
Further, the solution of the Cauchy problem (\ref{pr2}) is the distribution given by the formula
\begin{equation*}
x=e^{a\zeta_\tau}(ab\alpha\theta_\tau-b\delta_\tau^\alpha),
\end{equation*}
and the distribution
\begin{equation*}
x=e^{a\zeta_\tau}(\alpha a^2b\theta_\tau-2\alpha ab\delta_\tau^++b\delta_\tau^{\prime\alpha})
\end{equation*}
is the solution of the Cauchy problem (\ref{pr3}). Note that the solutions of the Cauchy problems (\ref{pr2}) and (\ref{pr3}) \textit{depend} on the value of $\alpha \in \mathbf C$.
\end{example}

\begin{example}
\label{example2}
Suppose that $I=\mathbf R$, and we are given a countable set $\Upsilon \subset (0,\infty)$. Let $\{a_\gamma\}_{\gamma \in \Upsilon}$, $\{b_\gamma\}_{\gamma \in \Upsilon} \subset \mathbf C$ be such that $\sum_{\gamma \in \Upsilon}|a_\gamma|$, $\sum_{\gamma \in \Upsilon}|b_\gamma|<\infty$. We define
\begin{equation*}
a=\sum_{\gamma \in \Upsilon} a_\gamma\theta_\gamma, \quad b=\sum_{\gamma \in \Upsilon} b_\gamma\theta_\gamma.
\end{equation*}
Let us consider the following Cauchy problem
\begin{equation}
\label{pr4}
x'=a(t)x+b'', \quad x=0 \text{ for }t<0,
\end{equation}
the second derivative of $b$ is specified by the formula
$b''=\sum_{\gamma \in \Upsilon} b_\gamma \delta_\gamma^{\prime\alpha(\gamma)}$, 
where $\alpha:I \mapsto \mathbf C$ is a bounded continuous function (Cauchy problem (\ref{pr2}) is a special case of Cauchy problem (\ref{pr4})). The solution of the Cauchy problem (\ref{pr4}) is the distribution given by the formula
\begin{equation*}
x=\sum_{\gamma \in \Upsilon} b_\gamma \exp\left( \int_{\gamma}^t a(s)ds \right)\left(\bigl(\alpha(\gamma)a(\gamma+)+\bigl(1-\alpha(\gamma)\bigr)a(\gamma-)\bigr)\theta_\gamma-\delta_\gamma^{\alpha(\gamma)} \right).
\end{equation*}
Observe that $x$ does not depend on values of $\alpha$ at the points of $\Upsilon$ where $b_\gamma \ne 0$ if and only if $a$ is continuous on $\Upsilon \cap \{\gamma: b_\gamma \ne 0\}$ (e.g., when $\Upsilon$ is the Cantor set, and $a$ is the Cantor function).
\end{example}

As follows from Theorem 3.5 below, the solutions of the Cauchy problems (\ref{pr1})--(\ref{pr4}) can also be obtained if the delta-functions and their derivatives in the-right hand sides of the differential equations in (\ref{pr1})--(\ref{pr4}) are replaced by converging families of locally-summable functions.

We also consider linear differential equations of higher orders, that is,
\begin{equation}
\label{ieqmm}
X^{(m)}-A_{m-1}X^{(m-1)}-\dots-A_0 X=F,
\end{equation}
where $A_i$ are the matrix-valued functions which are regulated (together with their derivatives of all orders defined almost everywhere) and $F$ is a matrix-valued distribution in $\mathcal R'$. We prove the well-posedness of the Cauchy problem for (\ref{ieqmm}) (Theorems 3.3'--3.5'), and so, is a sense, generalize the result in \cite{Der1} where the technique of quasidifferential equations is used to obtain the sufficient conditions on distribution $F \in \mathcal D'$ and coefficients $A_k$ for existence and uniqueness of the locally-summable solution of the Cauchy problem for (\ref{ieqmm}).

In what follows, we show that every distribution in $\mathcal D'$ admits a linear continuous extension from $\mathcal D$ to $\mathcal R$ (Theorem 2.8). Using Gelfand homomorphism induced by Banach algebra of regulated functions, we show  that the space of distributions $\mathcal R'$ is isomorphic to the space $\mathcal D'(I_*)$ of distributions with the smooth test functions defined on a totally-disconnected Hausdorff space (Theorem 2.19), so in this sense the use of the term "discontinuous test functions" is purely conventional
(nevertheless, the space of distributions $\mathcal D'(I_*)$ seems to be not studied before). \\[2mm]
\textbf{Acknowledgments.}
The author would like to thank Professor V.~Derr for his valuable comments, his support and his continuing interest in this work.

\section{Distributions}

Let us start with the definitions of functions spaces and functions algebras which are used throughout this paper.
We denote by $\mathbb L_{\loc}=\mathbb L_{\loc}(I)$ the linear space of locally-summable functions $I \mapsto \mathbf C$, where $I=(a,b) \subset \mathbf R$ is an open interval (in particular, $I=\mathbf R$). 
Let $\mathbb L_{\infty}=\mathbb L_{\infty}(I)$ be the algebra of functions essentially bounded on $I$, endowed with the norm
\begin{equation*}
\|g\|_{\mathbb L_{\infty}}=\esssup_{t \in I} |g(t)|.
\end{equation*}
In what follows, we denote by $\mathbb G \subset \mathbb L_{\infty}$ the algebra of \textit{regulated functions}, that is, the algebra of functions $g:I \mapsto \mathbf C$ possessing both one-sided limits
\begin{equation*}
g(t-):=\esslim_{s \to t-} g(s), \quad g(t+):=\esslim_{s \to t+} g(s)
\end{equation*}
for every $t \in I$ (equivalently, possessing at most first-kind discontinuities on $I$). As is well-known, $\mathbb G$ is a Banach algebra \cite{Dieu}.

\begin{theorem}[\cite{Der2}]
For every $g \in \mathbb G$ the set of points of discontinuity
$T(g):=\{t \in I: \sigma_t(g):=g(t+)-g(t-) \ne 0\}$
is at most countable.
\end{theorem}

We define support of the function $g \in \mathbb G$ to be the set
$$\supp(g)=\cl\{t \in I: g(t-) \ne 0 \text{ or } g(t+) \ne 0\},$$
where $\cl$ stands for the closure in $I$.

Further, let us denote by $\mathbb G^\infty \subset \mathbb L_\infty$ the subalgebra of functions $g \in \mathbb G$ such that for every $k \in \mathbf N$ there exists a regulated function $g^{(k)} \in \mathbb G$ (called the $k$-th derivative of $g$) such that
\begin{equation*}
g^{(k)}(t \pm)=\underset{s \to t \pm,~s \ne t}{\esslim} \left(\frac{g^{(k-1)}(s)-g^{(k-1)}(t\pm)}{s-t}\right), \quad t \in I,
\end{equation*}
where $g^{(0)}:=g$. As follows from the remark above, for every $g \in \mathbb G^\infty$, $k \in \mathbf N_0$, the set of points of discontinuity of the $k$-th derivative
$T(g^{(k)}):=\{t \in I: \sigma_t(g^{(k)}) \ne 0\}$
is at most countable.
We endow algebra $\mathbb G^\infty$ with the countable family of norms
\begin{equation*}
\|g\|_k=\max_{0 \leqslant i \leqslant k} \|g^{(i)}\|_{\mathbb L_\infty}, \quad k \in \mathbf N_0:=\mathbf N \cup \{0\}.
\end{equation*}

If we denote by $\mathbb C \subset \mathbb L_{\infty}$ the algebra of bounded continuous functions $I \mapsto \mathbf C$, then $\mathbb C$ is a proper subalgebra of $\mathbb G$. Furthermore, if $\mathbb C^\infty \subset \mathbb L_{\infty}$ is the algebra of
bounded infinitely (continuously) differentiable functions, then $\mathbb C^\infty$ is a proper subalgebra of $\mathbb G^\infty$ (we add "continuously" here since the elements of algebra $\mathbb G^\infty$ are also "infinitely differentiable" in the above sense).

In what follows, we denote by $\mathcal D=\mathcal D(I)$ the classical space of $\mathbb C^\infty$-test functions. Let $\mathcal D'=\mathcal D'(I)$ be the space of linear continuous functionals $\mathcal D \mapsto \mathbf C$ (called \textit{distributions}), endowed with weak* topology (see \cite{Gelf,Shi}). 

We proceed now to the definition of the space of $\mathbb G^\infty$-test functions, containing as a subspace the classical space of $\mathbb C^\infty$-test functions. Namely, let $\mathfrak J=\{J\}$ be the family of all subintervals in $\mathbf R$ such that $\bar{J} \subset I$. Every $\varphi \in \mathbb G^\infty(J)$, $J \in \mathfrak J$, can be extended to $I$ by assigning zero values on $I \setminus J$, so we may define
\begin{equation*}
\mathcal R=\bigcup_{J \in \mathfrak J} \mathbb G^\infty(J)
\end{equation*}
-- the linear space consisting of functions $\varphi \in \mathbb G^\infty$ possessing compact support in $I$. We endow $\mathcal R=\mathcal R(I)$ with the locally-convex topology of inductive limit of countably-normed spaces $\mathbb G^\infty(J)$, $J \in \mathfrak J$, and call its elements analogously as \textit{test functions}. 

The following theorem gives the description of the topology in $\mathcal R$ in terms of convergent sequences (its proof, as well as the proof of Lemma \ref{lem2} and the proof of Theorem \ref{thm3}, is provided in the last section).

\begin{lemma}
\label{lem12}
Given $\{\varphi_k\}_{k=1}^\infty \subset \mathcal R$ and $\varphi \in \mathcal R$, we have that $\varphi_k \to \varphi$ in $\mathcal R$ if and only if $\varphi_k \to \varphi$ in $\mathbb G^\infty$ and there exists $J \in \mathfrak J$ such that $\supp(\varphi_k) \subset J$ for all $k \in \mathbf N$. 
\end{lemma}

As one of the consequences of the characterization of topology in $\mathcal R$ given in Lemma \ref{lem12} we obtain that the classical space $\mathcal D$ of the smooth test functions is a subspace of $\mathcal R$.

Let $\mathcal R'=\mathcal R'(I)$ be the space of linear continuous functionals $\mathcal R \mapsto \mathbf C,$ endowed with weak* topology, whose elements are called analogously as \textit{distributions}. 
By definition of the inductive limit topology \cite{Rob}, 
a linear functional $f:\mathcal R \mapsto \mathbf C$ is a distribution if and only if given any $J \in \mathfrak J$ its restriction $f|_{\mathbb G^\infty(J)}$ is continuous.

As an example, given $f \in \mathbb L_{\loc}$, we may define distribution $f \in \mathcal R'$ whose value on the test function $\varphi \in \mathcal R$ is determined by the formula
\begin{equation*}
(f,\varphi)=\int_I f(t)\varphi(t)dt.
\end{equation*}
Note that the map $\mathbb L_{\loc} \mapsto \mathcal R'$ defined above is injective, since $\mathcal D$ is a subspace of $\mathcal R$, and analogous statement is true in $\mathcal D'$.

Let us consider some other examples of distributions in $\mathcal R'$.

\begin{example}[\cite{DK}]
\label{ex1}
Given $\tau \in I$, we define the \textit{right} and the \textit{left} delta-functions 
\begin{equation*}
(\delta_\tau^+,\varphi):=\varphi(\tau+), \quad (\delta_\tau^-,\varphi) := \varphi(\tau-),
\end{equation*}
where $\varphi \in \mathcal R$.
In general, given $\alpha \in \mathbf C$, we define
$$\delta_\tau^\alpha := \alpha\delta_\tau^++(1-\alpha)\delta_\tau^-.$$
Clearly, if $\varphi \in \mathcal D$, then $(\delta_\tau^\alpha,\varphi)=\varphi(\tau),$ so $\delta_\tau^\alpha \in \mathcal R'$ is the extension of the classical delta-function $\delta_\tau \in \mathcal D'$ from $\mathcal D$ to $\mathcal R$. 
Let us note that together with the family of delta-functions concentrated at $\tau \in I$ there exists the family of corresponding delta-sequences: if $\chi_S$ is the characteristic function of the set $S$, then
\begin{equation*}
f_k^\alpha=k\left(\alpha\chi_{(\tau,\tau+\frac{1}{2k})}+(1-\alpha)\chi_{(\tau-\frac{1}{2k},\tau)}\right) \to \delta_\tau^\alpha
\end{equation*}
in $\mathcal R'$.
Also, note that $f_k^1-f_k^0 \to \delta_\tau^+-\delta_\tau^-$ in $\mathcal R'$, while $f_k^1-f_k^0 \to 0$ in $\mathcal D'$.
\end{example}

\begin{example}
Let $k \in \mathbf N_0$, $\alpha \in \mathbf C$. Let us define the following distributions:
\begin{equation*}
(\delta_\tau^{(k)+},\varphi):= (-1)^k \varphi^{(k)}(\tau+), \quad (\delta_\tau^{(k)-},\varphi):= (-1)^k \varphi^{(k)}(\tau-), \quad \varphi \in \mathcal R,
\end{equation*}
and, in general,
\begin{equation}
\label{delderiv}
\delta_\tau^{(k)\alpha}:= \alpha \delta_\tau^{(k)+}+(1-\alpha)\delta_\tau^{(k)-},
\end{equation}
which we call the \textit{$k$-th derivatives of the delta-functions}. 
In fact, we haven't defined the derivative of a distribution yet. As it will be shown below, the distributions defined by (\ref{delderiv}) are indeed the derivatives of the delta-functions in $\mathcal R'$.
\end{example}

Suppose that $g \in \mathbb G^\infty$, $f \in \mathcal R'$. We define the product $gf \in \mathcal R'$ by the formula (see \cite{DK})
\begin{equation}
\label{product}
(gf,\varphi)=(fg,\varphi) := (f,g\varphi), \quad \varphi \in \mathcal R,
\end{equation}
where $g\varphi \in \mathcal R$. The operation of multiplication defined by (\ref{product}) is associative in the sense that the equality $(gh)f=g(hf)$ holds for any $g$, $h \in \mathbb G^\infty$, $f \in \mathcal R'$, and, clearly, coincides with the ordinary one for the regular distributions.
Furthermore, as follows from the next theorem, the operation of multiplication in $\mathcal R'$ is continuous.

\begin{theorem}
\label{thm3}
Given $g_k \to g$ in $\mathbb G^\infty$, $f_k \to f$ in $\mathcal R'$, we have $g_kf_k \to gf$ in $\mathcal R'$.
\end{theorem}

\begin{example} Let us find the product of the Heaviside function $\theta_\tau \in \mathbb G^\infty$ and the delta-function $\delta_\tau^\alpha \in \mathcal R'$. According to definition (\ref{product}) we have
\begin{equation*}
(\theta_\tau\delta_\tau^\alpha,\varphi)=(\delta_\tau^\alpha,\theta_\tau\varphi)=\alpha\theta_\tau(\tau+)\varphi(\tau+)+(1-\alpha)\theta_\tau(\tau-)\varphi(\tau-)=\alpha\varphi(\tau+),
\end{equation*}
where $\varphi \in \mathcal R$, so
\begin{equation*}
\theta_\tau\delta_\tau^\alpha=\alpha\delta_\tau^+.
\end{equation*}
Analogously, we have the following identities:
\begin{equation*} \theta_\tau\delta_\tau^{\prime\alpha}=\alpha\delta_\tau^{\prime+}, \quad \zeta_\tau\delta_\tau^\alpha=0, \quad \zeta_\tau\delta_\tau^{\prime\alpha}=-\alpha\delta_\tau^+.
\end{equation*}
\end{example}

\subsection{Properties of restriction and differentiation operators}

Let $$\Gamma: \mathcal R' \mapsto \mathcal D'$$ be the linear and, clearly, continuous operator of restriction from $\mathcal R$ to $\mathcal D$. 

\begin{theorem}
\label{extteo}
The operator $\Gamma$ is surjective, that is,
every distribution in $\mathcal D'$ has a linear continuous extension from $\mathcal D$ to $\mathcal R$. 
\end{theorem}
\begin{proof}
Indeed, since $\mathcal R$ is a locally-convex topological linear space, and $\mathcal D$ is a subspace of $\mathcal R$, the required extension exists by Hahn-Banach Theorem \cite{Kan}.
\end{proof}

As follows from the examples provided, such extension is always non-unique (e.g., $f$, $f+(\delta_\tau^+-\delta_\tau^-) \in \mathcal R'$ are two different extensions of the same distribution in $\mathcal D'$).
Denote $$\ker(\Gamma)=\{f \in \mathcal R': \Gamma(f)=0\}.$$ As follows from Example \ref{ex1}, $\ker(\Gamma) \setminus \{0\} \ne \varnothing$. Furthermore, $\ker(\Gamma)$ is a closed subspace of $\mathcal R'$.

In order to formulate the next theorem, which describes the structure of $\ker(\Gamma)$ and plays crucial role in our further considerations, we will need one definition and one supplementary statement.

First, given $f \in \mathcal R'$, we define the support $\supp(f) \subset I$ to be the minimal closed set such that for every $\varphi \in \mathcal R$ satisfying $\supp(\varphi) \cap \supp(f)=\varnothing$ we have $(f,\varphi)=0$. 

Second, given $k \in \mathbf N_0$, we define $\mathbb F_k \subset \mathbb G^\infty$ to be the subspace consisting of functions $g$ such that $g^{(i)} \in \mathbb C$ for all $i \ne k$, $g^{(k)}$ is piece-wise continuous. Let $\mathbb F \subset \mathbb G^\infty$ be the linear space spanned by $\mathbb F_k$, $k \in \mathbf N_0$. Let $\mathcal R_F := \mathcal R \cap \mathbb F$ be endowed with the topology induced by the topology of $\mathcal R$, so that $\mathcal R_F$ is a subspace of $\mathcal R$.
The following statement is technical but essential.

\begin{lemma}
\label{lem2}
The subspace $\mathcal R_F$ is dense in $\mathcal R$. Furthermore, for every $\tau \in I$, $\varphi \in \mathcal R$ there exists a sequence $\{\varphi_l\}_{l=1}^\infty \subset \mathcal R_F$ such that $\varphi_l \to \varphi$
in $\mathcal R$, $\varphi_l^{(j+1)}$ is continuous in $t=\tau$ for all $j \geqslant l$, $l \in \mathbf N$, and $$|\varphi^{(j)}_l(\tau \pm )-\varphi^{(j)}(\tau\pm)|<l^{-1},$$ where $l \geqslant j$, $j \in \mathbf N_0$
\end{lemma}

\begin{theorem} \label{teo1} 1) Let $f \in \ker(\Gamma)$. Then there exist uniquely determined functions $t \mapsto c^k_f(t)$, $k \in \mathbf N_0$, such that for every $\varphi \in \mathcal R$
\begin{equation}
\label{rep3}
(f,\varphi)=\sum_{k \in \mathbf N_0} \sum_{t \in I}c^k_f(t) \sigma_t\bigl(\varphi^{(k)}\bigr)
\end{equation}
where the distribution
\begin{equation}
\label{rep1}
\mathcal R \ni \varphi \mapsto \sum_{k \in \mathbf N_0} \sum_{t \in I}c^k_f(t) \sigma_t\bigl(\varphi^{(k)}\bigr)
\end{equation}
is defined to be the extension of the functional
\begin{equation}
\label{rep2}
\mathcal R_F \ni \varphi \mapsto  \sum_{k \in \mathbf N_0}\sum_{t \in I}c^k_f(t) \sigma_t\bigl(\varphi^{(k)}\bigr)
\end{equation}
from $\mathcal R_F$ to $\mathcal R$. 

2) If $f \in \ker({\Gamma})$, $t \in I$, then $c_f^k(t)=0$ starting with certain $k$.

3) If $f \in \ker({\Gamma})$, $t \in I \setminus \supp(f)$, then $c_f^k(t)=0$ for all $k \in \mathbf N_0$.
\end{theorem}

\begin{remark}
\label{rem1}
Observe that every $\varphi \in \mathcal R_F$ has continuous derivatives starting with certain $K \in \mathbf N_0$, and the sets of points of discontinuity $T(\varphi^{(k)})$ ($0 \leqslant k \leqslant K$) are finite, so in (\ref{rep2}) there is only finite number of non-zero summands. Also, note that the extension of (\ref{rep2}) exists (take $f$) and unique, since $\mathcal R_F$ is dense in $\mathcal R$ (Lemma \ref{lem2}). 
\end{remark}

\begin{proof}
1) Let us show that the equality (\ref{rep3}) holds for every $\varphi \in \mathcal R_F$. 
Let $\mathcal R_F(k,\tau)$ be the subspace of $\mathcal R_F$ consisting of all test functions $\varphi$ such that $\varphi^{(i)} \in \mathbb C$ for all $i \in \mathbf N_0$, $i \ne k$, and $\varphi^{(k)}$ may have discontinuity only at $\tau \in I$. 
Let us show that for any $k \in \mathbf N_0$ and $\tau \in I$ we may find $a \in \mathbf C$ such that
\begin{equation*}
(f,\varphi)=a\sigma_\tau(\varphi^{(k)}), \quad \varphi \in \mathcal R_F(k,\tau).
\end{equation*}
Suppose that $\varphi_0 \in \mathcal R_F(k,\tau)$ is such that $\sigma_\tau(\varphi_0^{(k)})=1$. We define: $$a := (f,\varphi_0).$$ Let us show that the value of $a$ does not depend on choice of $\varphi_0$. Indeed, suppose that $\varphi_1 \in \mathcal R_F(k,\tau)$, $\sigma_\tau(\varphi_1^{(k)})=1$; then $\varphi_0-\varphi_1 \in \mathcal D$ and, since $f\in \ker(\Gamma)$, we obtain that $(f,\varphi_0-\varphi_1)=0$, i.e., $(f,\varphi_1)=a$. 
Now let $\varphi \in \mathcal R_F(k,\tau)$ be arbitrary. Then $\varphi=\sigma_\tau(\varphi)\varphi_2$ for certain $\varphi_2 \in \mathcal R_F(k,\tau)$ such that $\sigma_\tau(\varphi_2^{(k)})=1$, so due to linearity of $f$ we have
\begin{equation*}
(f,\varphi)=\sigma_\tau(\varphi)(f,\varphi_2)=a\sigma_\tau(\varphi).
\end{equation*}
Let us define: $$c_f^k(t) := a.$$ 
Suppose that $\varphi \in \mathcal R_F$. According to Remark \ref{rem1} there exist a number $K \in \mathbf N_0$ and the functions $\varphi_{ik} \in \mathcal R_F(k,\tau_{ik})$, where $\{\tau_{ik}\}_{i=1}^{m_k}=T(\varphi^{(k)})$, $0 \leqslant k \leqslant K$, such that
\begin{equation*}
\varphi-\sum_{k=0}^K \sum_{i=1}^{m_k}\varphi_{ik} \in \mathcal D,
\end{equation*}
so $\sigma_{\tau_i}(\varphi^{(k)})=\sigma_{\tau_i}(\varphi_{ik})$, $1 \leqslant i \leqslant m_k$, $0 \leqslant k \leqslant K$. Consequently, we have
\begin{equation*}
(f,\varphi)=\sum_{k=0}^K \sum_{i=1}^{m_k}(f,\varphi_{ik})=\sum_{k=0}^K \sum_{i=1}^{m_k}c_f^k(\tau_i)\sigma_{\tau_i}(\varphi^{(k)})=\sum_{k \in \mathbf N_0} \sum_{t \in I}c_f^k(t)\sigma_{t}(\varphi^{(k)}),
\end{equation*}
where the latter equality is obtained by adding the zero summands.

2) Suppose the contrary. Given $t \in I$, we denote $a_k=c_f^k(t)$. First, suppose that the sequence $\{a_k\}_{k=1}^\infty$ is bounded. Without loss of generality we may assume that $a_k \ne 0$, $k \in \mathbf N_0$. According to \cite{Shi} there is the test function $\varphi \in \mathcal D$ such that
\begin{equation*}
|\varphi^{(k)}(t)|>|a_k|^{-1}, \quad k \in \mathbf N_0.
\end{equation*}
Let $\theta_t$ be the Heaviside function which is discontinuous at $t$. Clearly, $\theta_t \in \mathbb G^\infty$. Then $\theta \varphi \in \mathcal R$, and we have
\begin{equation*}
|(\theta_t\varphi)^{(k)}(t+)|>|a_k|^{-1}, \quad (\theta_t\varphi)^{(k)}(t-)=0, \quad k \in \mathbf N_0.
\end{equation*}
According to Lemma \ref{lem2} there exists the sequence $\{\psi_l\}_{l=1}^\infty \subset \mathcal R$ such that $\psi_l^{(j+1)}$ is continuous at $t$ for all $j \geqslant l$, $\psi_l \to \theta\varphi$
in $\mathcal R$, and 
\begin{equation*}
|\psi^{(j)}_l(t+)-(\theta_t\varphi)^{(j)}(t+)|<l^{-1}, \quad |\psi^{(j)}_l(t-)-(\theta_t\varphi)^{(j)}(t-)|<l^{-1}
\end{equation*}
for all $l \geqslant j$, $j \in \mathbf N_0$, so
\begin{equation*}
|\psi^{(j)}_l(t+)|>|a_l|^{-1}-l^{-1}, \quad |\psi^{(j)}_l(t-)|<l^{-1}, 
\end{equation*}
for all $l \geqslant j$, $j \in \mathbf N_0$. Without loss of generality we may assume that $\psi^{(j)}_l$ is continuous on $I \setminus \{t\}$, $j \in \mathbf N_0$. Due to continuity of $f$ there exists the limit $$\lim(f,\psi_l)=(f,\theta\varphi).$$ Along with that, we have the equality
\begin{equation*}
(f,\psi_l)=\sum_{j=0}^{l} a_j \sigma_t(\psi^{(j)}_l),
\end{equation*}
and
\begin{equation*}
|(f,\psi_{l+1})-(f,\varphi_l)|=\bigl|a_{l+1}\sigma_t\bigl(\psi^{(l+1)}_{l+1}\bigr)\bigr|>|a_{l+1}|\left(|a_{l+1}|^{-1}-2l^{-1}\right)=1-2|a_{l+1}|l^{-1} \to 1,
\end{equation*}
so the limit $\lim (f,\psi_l)$ does not exists, which contradicts to our assumption. 

Second, suppose that the sequence $\{a_k\}_{k=1}^\infty$ is unbounded. Then we may choose the test function $\varphi \in \mathcal D$ such that
$|\varphi^{(k)}(t)|>1$, $k \in \mathbf N_0$,
and repeat the previous argument. Then
\begin{equation*}
|(f,\psi_{l+1})-(f,\varphi_l)|=\bigl|a_{l+1}\sigma_t\bigl(\psi^{(l+1)}_{l+1}\bigr)\bigr|>|a_{l+1}|\left(1-2l^{-1}\right),
\end{equation*}
so the limit $\lim (f,\psi_l)$ does not exist, which again leads to contradiction.

3) Since $\supp(f)$ is a closed set, given any $k \in \mathbf N_0$ there exists the test function $\varphi \in \mathcal R_F(k,t)$ such that $\sigma_t(\varphi^{(k)})=1$ and $\supp(\varphi) \cap \supp(f)=\varnothing$. Then
\begin{equation*}
c_f^k(t)=(f,\varphi),
\end{equation*}
but $(f,\varphi)=0$, so $c_f^k(t)=0$.
\end{proof}

\begin{example}
Suppose that in (\ref{rep3}) $c^k_f(t)=0$ for $k>0$, $t \in I$ and for $k=0$, $t \ne \tau$, while $c^0_f(\tau)=a$. Then
\begin{equation}
\label{f1}
f=a(\delta_\tau^+-\delta_\tau^-).
\end{equation}
Suppose that in (\ref{rep3}) $c^k_f(t)=0$ if $k \ne m$, $t \in I$, and $c^m_f(t)=1$ if $t \in I$, $m \in \mathbf N_0$. Then
\begin{equation}
\label{f2}
(f,\varphi)=\int_I \varphi^{(m+1)}(t)dt, \quad \varphi \in \mathcal R.
\end{equation}
Clearly, both distributions (\ref{f1}), (\ref{f2}) are in $\ker(\Gamma)$.
\end{example}

Let $\Delta \subset \ker(\Gamma)$ be the subspace consisting of the test functions of the form
\begin{equation*}
\varphi \mapsto \sum_{t \in I} c(t)\sigma_t(\varphi), \quad \varphi \in \mathcal R,
\end{equation*}
where $I \ni t \mapsto c(t) \in \mathbf C$ (equivalently, $\Delta$ consists of the elements $f \in \ker(\Gamma)$ such that $c^k_f \equiv 0$ if $k \geqslant 1$).
Suppose that $\{f_k\}$ is a sequence of locally absolutely continuous functions $I \mapsto \mathbf C$ such that
\begin{equation*}
f_k \to 0, \quad f'_k \to g
\end{equation*}
in $\mathcal R'$. Then it is natural to consider distribution $g \in \mathcal R'$ as the derivative of the zero distribution in $\mathcal R'$. Let us show that $g \in \Delta$. Observe first that $g \in \ker(\Gamma) \supset \Delta$, since $f'_k \to 0$ in $\mathcal D'$. Further, let $\varphi \in \mathcal R_F(0,\tau)$, $\tau \in I$. Then
\begin{multline}
\notag
(g,\varphi)=\lim \int_I f'_k(t)\varphi(t)dt=\lim\left(\int_{a}^\tau f'_k(t)\varphi(t)dt+\int_{\tau}^b f'_k(t)\varphi(t)dt \right)=\\
=\lim\left(-f_k(\tau)\sigma_\tau(\varphi)- \int_I f_k(t)\varphi'(t) \right)=-\sigma_\tau(\varphi)\lim f_k(\tau),
\end{multline}
where $\lim f_k(\tau)$ exists since the value $(g,\varphi)$ is defined. Let $c(t):=- \lim f_k(\tau)$. 
Suppose that we are given $\varphi \in \mathcal R_F(j,\tau)$, $j \geqslant 1$. Then
\begin{equation*}
(g,\varphi)=\lim \int_I f'_k(t)\varphi(t)dt=-\lim \int_I f_k(t)\varphi'(t)dt=0
\end{equation*}
for every $\tau \in I$. Consequently, if $\varphi \in \mathcal R_F$, then
\begin{equation}
\label{gen}
(g,\varphi)=\sum_{t \in T(\varphi)} c(t)\sigma_t(\varphi)=\sum_{t \in I} c(t)\sigma_t(\varphi),
\end{equation}
where the latter equality is obtained by adding zero summands. The distribution (\ref{gen}) can be extended uniquely to $\mathcal R$, by the definition
$g \in \Delta$.

Conversely, let us show that given a distribution $g \in \Delta$ and a sequence of locally summable functions $\{g_k\}_{k=1}^\infty$ such that $g_k \to g$ in $\mathcal R'$, in assumption that the sequence of functions (i.e., primitives) $\{f_k\}_{k=1}^\infty$ defined by the formula $$t \mapsto f_k(t):=\int_{t_0}^t g_k(s)ds \quad (t_0 \in I)$$ tends to a regular distribution $f \in \mathcal R'$, we have that $f$, as an ordinary function, is identically equal to a constant. 
If we show that there exists $r \in \mathbf C$ such that 
$$(f_k,\varphi) \to r\int_I \varphi(t)dt,$$ 
for every $\varphi \in \mathcal R$, then 
the proof would be complete.
Indeed, given $\varphi \in \mathcal R$ (suppose that $\supp(\varphi) \subset [u,v] \subset I$), we have
\begin{multline}
\notag
\int_I f_k(t)\varphi(t)dt=\int_{u}^v \left(\int_{t_0}^t g_k(s)ds\right)d\int_u^t \varphi(s)ds=\\=\left(\int_{t_0}^v g_k(s)ds \right) \int_u^v \varphi(s)ds-\int_u^v \left(\int_u^t \varphi(s)ds \right) g_k(t)dt \to \\
\to \bigl(g,\chi_{(t_0,v)}\bigr)\int_u^v \varphi(s)ds-\left(g,\int_u^t\varphi(s)ds\chi_{(u,v)}(t) \right)=\\=\bigl(c_g(v)-c_g(t_0)\bigr)\int_I \varphi(s)ds-c_g(v)\int_I \varphi(s)ds=-c_g(t_0)\int_I \varphi(s)ds,
\end{multline}
so we may put $r:=-c_g(t_0)$. The proof is complete.

Taking into account the above considerations, we define the derivative of the distribution $f \in \mathcal R'$ to be any distribution $f' \in \mathcal R'$ such that
\begin{equation}
\label{deriv}
f' \in \frac{df}{dt}+\Delta, \text{ where }\left( \frac{df}{dt},\varphi \right) := -(f,\varphi'), \quad \varphi \in \mathcal R,
\end{equation}
where $\varphi'$ is assumed to be defined almost everywhere.
We denote by $D(f) \subset \mathcal R'$ the family of all derivatives of the distribution $f \in \mathcal R'$.
As it immediately follows from the definition, the operator of differentiation $f \mapsto D(f)$ is multi-valued.

We define the derivatives of higher orders inductively.

\begin{example}
Let $\theta_\tau$ be Heaviside function. Then
\begin{equation*}
\left(\frac{d\theta_\tau}{dt},\varphi\right)=-\int_\tau^b \varphi'(t)dt=\varphi(\tau+),
\end{equation*}
so $(\frac{d\theta_\tau}{dt},\varphi)=\delta_\tau^+$. Consequently, for every $\alpha \in \mathbf C$ delta-function $\delta_\tau^\alpha$ is the derivative of $\theta_\tau$.
\end{example}

\begin{example}
As follows from definition (\ref{deriv}) the distribution (\ref{delderiv}) is indeed the $k$-th derivative of $\delta_\tau^\alpha$. Furthermore,
\begin{equation*}
\delta_\tau^{\prime\alpha} + \beta(\delta_\xi^+-\delta_\xi^-)
\end{equation*}
is the derivative of $\delta_\tau^\alpha$ for any $\beta \in \mathbf C$, $\xi \in I$,
\begin{equation*}
\delta_\tau^{\prime\alpha} + \beta(\delta_\xi^{\prime+}-\delta_\xi^{\prime-})
\end{equation*}
is the derivative of $\delta_\tau^\alpha+\beta(\delta_\xi^+-\delta_\xi^-)$, and  $$\delta_\tau^\alpha+\beta(\delta_\xi^+-\delta_\xi^-)$$ is the derivative of $\theta_\tau$ for any $\beta \in \mathbf C$, $\xi \in I$.
\end{example}

\begin{example}
\label{impex}
Let $f:I \mapsto \mathbf C$ be a locally absolutely continuous function. Then
\begin{multline}
\label{dotf}
(f',\varphi)=\int_I f'(t)\varphi(t)dt=\int_I \varphi(t)d f(t)=\\=-\int_I f(t)d\varphi(t)=-\int_I f(t)\varphi'(t)dt-\sum_{t \in I} f(t)\sigma_t(\varphi), \quad \varphi \in \mathcal R
\end{multline}
where, as usual, the latter sum is defined to be the extension of the corresponding functional from $\mathcal R_F$ to $\mathcal R$. Since $\sum_{t \in I} f(t)\sigma_t(\varphi) \in \Delta$, we have
\begin{equation*}
f'+\Delta=D(f),
\end{equation*}
where $f'$ is defined by (\ref{dotf}).
\end{example}

\begin{theorem}
\label{importantremark}
Suppose that $f \in \mathcal R'$, and we are given its $m$-th derivative $f^{(m)} \in \mathcal R'$. Then the intermediate derivatives $f^{(k)}$, $0 \leqslant k \leqslant m-1$, are determined uniquely. 
\end{theorem}
\begin{proof}
We will prove this result by induction over $k$. By definition $f^{(0)}=f$ is determined uniquely. Suppose that $f^{(i)}$ are uniquely determined for $0 \leqslant i \leqslant k-1$. Then, since the primitive of a distribution is defined uniquely up to a constant summand, we have that $f^{(k)}$ is defined uniquely by $f^{(m)}$ up to a polynomial of degree $m-k-1$. Along with that, if $f^{(k)}_1$, $f^{(k)}_2$ are two different derivatives of $f^{(k-1)}$, then $$f^{(k)}_1-f^{(k)}_2 \in \Delta.$$ Now since $\Delta$ does not contain polynomials except the one identically equal to zero, we have that $f^{(k)}$ is uniquely determined.
\end{proof}

\subsection{Theorems on structure of distributions} The following statement is a generalization of well-known result on structure of distributions in $\mathcal D'$ \cite{Shi}.

\begin{theorem}
\label{struc}
Let $f \in \mathcal R'$. Then there exist functions $f_k \in \mathbb L_{\loc}$ and functions $c^k:I \mapsto \mathbf C$, $k \in \mathbf N_0$, such that
\begin{equation}
\label{struceq}
(f,\varphi)=\sum_{k \in \mathbf N_0} (-1)^k\int_I f_k(t)\varphi^{(k)}(t)dt+\sum_{k \in \mathbf N_0} \sum_{t \in I} c^k(t)\sigma_t(\varphi^{(k)}), \quad \varphi \in \mathcal R.
\end{equation}
If $f=0$ for $t<t_0$, that is, for every test function $\varphi \in \mathcal R$ such that $\supp(\varphi) \subset (a,t_0)$ we have $(f,\varphi)=0$, then there exist $f_k$ such that $f_k=0$ for $t<t_0$.
\end{theorem}
\begin{proof}
According to \cite{Shi} there exist functions $f_k \in \mathbb L_{\loc}$, $k \in \mathbf N_0$, such that
\begin{equation}
\label{f}
(f,\varphi)=\sum_{k \in \mathbf N_0} (-1)^k\int_I f_k(t)\varphi^{(k)}(t)dt, \quad \varphi \in \mathcal D.
\end{equation}
(if $f=0$ for $t<t_0$, then there exist $f_k$ such that $f_k=0$ for $t<t_0$). If we consider in (\ref{f}) the test functions $\varphi \in \mathcal R$,
then we get certain extension of $\Gamma(f)$ from $\mathcal D$ to $\mathcal R$. According to Theorem \ref{teo1} the family of all extensions of $\Gamma(f)$ from $\mathcal D$ to $\mathcal R$ consists of all distributions of form (\ref{struceq}). Since $f$ is one of such extensions, there exist locally-summable functions $c^k$ such that (\ref{struceq}) is true.
\end{proof}

\begin{theorem}
\label{deltateo}
Let $f \in \mathcal R'$, $\supp(f)\subset \{\tau\}$. Then
\begin{equation}
\label{deltarepr}
f=\sum_{k=0}^K a_k \delta_\tau^{(k)\alpha_k}+\sum_{k=0}^M b_k \left(\delta_\tau^{(k)+}-\delta_\tau^{(k)-}\right)
\end{equation}
for certain $K$, $M \in \mathbf N_0$, $a_k$, $b_k \in \mathbf C$, $\alpha_k \in \mathbf C$.
\end{theorem}

As it follows from Theorem \ref{deltateo}, in $\mathcal R'$ there are no other extensions of delta-function and its derivatives from $\mathcal D$ to $\mathcal R$, concentrated at a single point, except defined in the examples above.

\begin{proof}
If $\supp(f)=\varnothing$, then $f=0$, and the proof is complete. So, we may suppose that $\supp(f)=\{\tau\}$.
As follows from the definition of support, we have the inclusion $\supp\bigl(\Gamma(f)\bigr) \subset \{\tau\}$. According to \cite{Shi} there exist $K \in \mathbf N_0$, $a_k \in \mathbf C$, $0 \leqslant k \leqslant K$, such that
\begin{equation*}
\Gamma(f)=\sum_{k=0}^K a_k\delta_\tau^{(k)}.
\end{equation*}
Let $\alpha_k \in \mathbf C$, $0 \leqslant k \leqslant K$. Then
\begin{equation*}
\varphi \mapsto \sum_{k=0}^K a_k\delta_\tau^{(k)\alpha_k}
\end{equation*}
is an extension of $\Gamma(f)$ from $\mathcal D$ to $\mathcal R$.
According to Theorem \ref{teo1} any extension of $\Gamma(f)$, concentrated at $\tau$, has form (\ref{deltarepr}) for certain $M \in \mathbf N_0$, $a_k$, $b_k \in \mathbf C$, $\alpha_k \in \mathbf C$. Since $f$ is one of such extensions, we obtain the statement of the theorem.
\end{proof}

Let us show that the space of distributions $\mathcal R'=\mathcal R'(I)$  is isomorphic (as a topological linear space) to the space of distributions with the smooth test functions defined on a certain totally disconnected set. Namely, let $I_*$ be the maximal ideal space of the Banach algebra $\mathbb G=\mathbb G(I)$, that is, the space of all continuous algebra homomorphisms $\mathbb G \mapsto \mathbf C$, endowed with weak* topology \cite{Gam}. As is shown in \cite{BK}, $I_*$ is a totally disconnected Hausdorff space which consists of homomorphisms of the form
$$g \mapsto g(t-), \quad g \mapsto g(t+),$$ where $g \in \mathbb G$, $t \in I$. We denote these homomorphisms by $t-$ and $t+$, respectively, and use notations $g(t-)$ for $(t-)(g)$ and $g(t+)$ for $(t+)(g)$, $g \in \mathbb G$. So, as the set, $I_*$ is in one-to-one correspondence with $I \times \{-1,1\}$, where $-1$ stands for the left hand-side limit evaluation homomorphism, and $1$ stands for the right hand-side limit evaluation homomorphism.
In what follows, we put $[t+]=[t-]:=t$.
Since algebra $\mathbb G$ is regular and symmetric, i.e.,
\begin{equation*}
\|g^2\|_{\mathbb L_\infty}=\|g\|_{\mathbb L_\infty}^2
\end{equation*}
and $\bar{g} \in \mathbb G$ for every $g \in \mathbb G$, we have that $\mathbb G$ is isometrically isomorphic to $\mathbb C(I_*)$ \cite{Gam}. Let $\mathbb C^\infty(I_*)$ be the algebra consisting of infinitely continuously differentiable functions $g \in \mathbb C(I_*)$, that is, the functions such that for every $k \in \mathbf N$, $t\cdot \in I_*$ there exists the limit
\begin{equation*}
g^{(k)}(t\cdot)=\underset{r \to t\cdot,~r \ne t\cdot}{\lim} \left(\frac{g^{(k-1)}(r)-g^{(k-1)}(t\cdot)}{[r]-t} \right)
\end{equation*}
and $g^{(k)} \in \mathbb C(I_*)$ (by definition, $g^{(0)}:=g$). 
We introduce in $\mathbb C^\infty(I_*)$ the countable family of norms
\begin{equation*}
\|g\|_k=\max_{0 \leqslant i \leqslant k} \sup_{t\cdot \in I_*}|g^{(i)}(t\cdot)|, \quad g \in \mathbb C^\infty(I_*), \quad k \in \mathbf N_0.
\end{equation*}

\begin{lemma}
$\mathbb G^\infty(I)$ is isometrically isomorphic to $\mathbb C^\infty(I_*)$.
\end{lemma}
\begin{proof}
Suppose that $f \in \mathbb G(I)$. In what follows, we denote the image of $f$  in $\mathbb C^\infty(I_*)$ under the map defined above by $f_*$. 
Let $g \in \mathbb G^\infty(I)$. Show that for every $k \in \mathbf N_0$ 
\begin{equation}
\label{reqeq}
(g^{(k)})_*=(g_*)^{(k)}. 
\end{equation}
Let us prove this statement inductively by $k$. Clearly, $$(g^{(0)})_*=(g_*)^{(0)}.$$ Suppose that $m \in \mathbf N$, and (\ref{reqeq}) is true for $k=m-1$. We have to show now that (\ref{reqeq}) is true for $k=m$. Observe that due to the fact that the set of points of discontinuity $T(g^{(m)})$ has zero measure, we have the equality
\begin{equation*}
g^{(m)}(t\pm)=\underset{s \to t \pm,~s \ne t}{\esslim} \left(\frac{g^{(m-1)}(s)-g^{(m-1)}(t\pm)}{s-t}\right)=\underset{s \to t \pm,~s \ne t}{\lim} \left(\frac{g^{(m-1)}(s\cdot)-g^{(m-1)}(t\pm)}{s-t}\right)
\end{equation*}
As it follows from our assumption and from the definition of the topology in $I_*$, the latter is equal to 
\begin{equation*}
\underset{r \to t \pm,~r \ne t}{\lim} \left(\frac{g_*^{(m-1)}(r)-g_*^{(m-1)}(t\pm)}{[r]-t}\right)=g_*^{(m)}(t \pm)
\end{equation*}
for every $t \in I$, so the equality (\ref{reqeq}) holds for $k=m$, so $\mathbb G^\infty(I)$ is isomorphic to $\mathbb C^\infty(I_*)$. Further, as follows from (\ref{reqeq}), we have that
\begin{equation*}
\|g\|_k=\|g_*\|_k, \quad k \in \mathbf N_0,
\end{equation*}
so $\mathbb G^\infty(I)$ is isometrically isomorphic to $\mathbb C^\infty(I_*)$.
\end{proof}

Let $\mathcal D(I_*)$ be the space consisting of elements $\varphi \in \mathbb C^\infty(I_*)$ such that there exist $c_\varphi$, $d_\varphi \in I$ possessing the property $\varphi(t \cdot )=0$ if $t<c_\varphi$, $\varphi(t\cdot)=0$ if $t>d_\varphi$, endowed with the convergence: $\varphi_k \to \varphi$ in $\mathcal D(I_*)$ if and only if $\varphi_k \to \varphi$ in $\mathbb C^\infty(I_*)$ and there exist $c$, $d \in I$ such that $\varphi_k(t \cdot )=0$ if $t<c_\varphi$, $\varphi_k(t\cdot)=0$ if $t>d_\varphi$ for all $k \in \mathbf N$. Let $\mathcal D'(I_*)$ be the space of all continuous linear functionals $\mathcal D(I_*) \mapsto \mathbf C$. As a simple corollary of the results above, we obtain the following statement.

\begin{theorem}
\label{isothm}
$\mathcal R'(I)$ is isomorphic to $\mathcal D'(I_*)$.
\end{theorem}

\section{Linear differential equations}

Let $\mathcal R^{n\prime}$ be the space of vector-valued distributions $\mathcal R \mapsto \mathbf C^n$ with the operations and convergence defined componentwise. We introduce analogous notations for the spaces of vector-valued and matrix-valued distributions and functions $I \mapsto \mathbf C^n$ and $I \mapsto \mathbf C^{n \times m}$, respectively.

Let us consider in the space $\mathcal R^{n\prime}$ the following linear differential equation
\begin{equation}
\label{e1}
x'-A(t)x=0
\end{equation}
with the matrix of coefficients $A \in \mathbb G^\infty_{n \times n}$.
Analogously, the solution of the differential equation (\ref{e1}) is the distribution $x \in \mathcal R^{n\prime}$ which possesses a derivative $x' \in \mathcal R^{n\prime}$ such that after the substitution of $x$ and $x'$ in (\ref{e1}) we obtain the identity $(x',\varphi)=(A(t)x,\varphi)$ for all $\varphi \in \mathcal R$.

\begin{theorem}
\label{teo2}
There are no other solutions in the space $\mathcal R^{n\prime}$ of the differential equation {\rm(\ref{e1})} except the classical ones. 
\end{theorem}
\begin{proof}
Without loss of generality we conduct our proof for the case $n=1$.
First, let us show that there are no other solutions of equation (\ref{e1}) being considered on the space of smooth test functions $\mathcal D$ except the classical ones. Let $x \in \mathcal R'$ be a solution of equation (\ref{e1}). By definition, the derivatives of $x$ coincide on $\mathcal D \subset \mathcal R$, so $(x',\varphi)=-(x,\varphi')$. Also note that if we are given a distribution $h \in \mathcal R'$ and a function $g \in \mathbb G^\infty$ which is continuous (and, thus, differentiable), then we have $((gh)',\varphi)=-(gh,\varphi')=-(h,g\varphi')=(h,-(g\varphi)'+g'\varphi)$ for all $\varphi \in \mathcal D$, that is,
we have the identity $(gh)'=h'g+hg'$ on $\mathcal D$. Now let us represent $x$ in the form $x=e^By$, where $B \in \mathbb G^\infty$ is a (continuous) primitive of $A$, and $y \in \mathcal R'$. Then
\begin{equation*}
x'=(e^B)'y+e^By',
\end{equation*}
on $\mathcal D$, so our equation is equivalent to
\begin{equation*}
Ae^By+e^By'=Ae^By,
\end{equation*}
which is, in turn, implies that $(e^By',\varphi)=0$ for all $\varphi \in \mathcal D$. We may multiply both parts of this equality by $e^{-By}$ to get the identity $(y',\varphi)=0$ for all $\varphi \in \mathcal D$. According to \cite{Shi} $y \equiv \const$. Thus, $x$ restricted to $\mathcal D$ can be only an ordinary function (that is, regular distribution).

We obtained that if $x$ is the solution of (\ref{e1}) in $\mathcal R'$, then 
$$x=x_0+y,$$
where $x_0$ is the classical solution, $y \in \ker(\Gamma)$. Clearly, $y$ is also a solution of equation (\ref{e1}), i.e.,
\begin{equation*}
\sum_{k \in \mathbf N_0}\sum_{t \in I} c_y^k(t)\sigma_t(\varphi^{(k+1)})+(z,\varphi)=\sum_{k \in \mathbf N_0}\sum_{t \in I} c_y^k(t)\sigma_t\left((A\varphi)^{(k)}\right), \quad \varphi \in \mathcal R,
\end{equation*}
where $z \in \Delta$. Let us choose certain $\tau \in I$, $\varphi \in \mathcal R(k,\tau)$. Then, since $c_y^k(\tau)=0$ starting with certain $K$ according to Theorem \ref{teo1}, we have the equality
\begin{equation}
\label{equ1}
\sum_{k=0}^K c_y^k(\tau)\sigma_\tau(\varphi^{(k+1)})+c_z(\tau)\sigma_\tau(\varphi)=\sum_{k=0}^K c_y^k(\tau)\sigma_\tau\left((A\varphi)^{(k)}\right).
\end{equation}
Since $\sigma(\varphi^{(K+1)})$ is not contained in the right-hand side of (\ref{equ1}), we find that $c_y^K(\tau)=0$. Consequently, the right-hand side of (\ref{equ1}) does not contain the summand corresponding to $k=K$. Then $\sigma_\tau(\varphi^{(K)})$ is not contained in the right-hand side of (\ref{equ1}), so $c_y^{K-1}(\tau)=0$. We may continue this process to obtain that $c_y^1(\tau)=0$. Then the right-hand side of (\ref{equ1}) does not contain $\sigma_\tau(\varphi')$ and, as a result, we get $c_y^0(\tau)=0$ (and also $c_z(\tau)=0$). Now since $\tau \in I$ was chosen arbitrarily, we obtain that $c_y^k(t)=0$ for all $t \in I$, $k \in \mathbf N_0$. Consequently, the only solution of differential equation (\ref{e1}) in $\ker(\Gamma)$ is the zero solution, so $x=x_0$.
\end{proof}

Further, let us consider in the space $\mathcal R'$ the differential equation
\begin{equation}
\label{b1}
x'=f,
\end{equation}
where $f \in \mathcal R'$. In what follows, we assume that $f=0$ if $t<t_0$. 
Analogously, we call the solution of the differential equation (\ref{b1}) the distribution $x \in \mathcal R'$ which possesses a derivative $x' \in \mathcal R'$ such that after substitution of $x'$ in (\ref{b1}) the equation (\ref{b1}) becomes the identity $(x',\varphi)=(f,\varphi)$, $\varphi \in \mathcal R$.

\begin{theorem}
\label{primteo}
There exists solution $x$ of equation {\rm(\ref{b1})} such that $x=0$ if $t<t_0$.
\end{theorem}

The solution $x$ of equation (\ref{b1}) is called the \textit{primitive} of distribution $f$.

\begin{proof}
According to Theorem \ref{struc}
there exist the locally-summable functions $f_k$ which are equal to 0 for $t<t_0$, and functions $c^k$, $k \in \mathbf N_0$ such that (\ref{struceq}) holds. Let us define (let $t_1<t_0$, $t_1 \in I$):
\begin{multline}
\label{xdef}
(x,\varphi):=-\int_{I}\left(\int_{t_1}^t f_0(s)ds\right) \varphi'(t)dt+\\+\sum_{k \in \mathbf N_0} (-1)^{k} \int_I f_{k+1}(t)\varphi^{(k)}(t)dt-\sum_{k \in \mathbf N_0}\sum_{t \in I} c^{k+1}(t)\sigma_t(\varphi^{(k)}), \quad \varphi \in \mathcal R.
\end{multline}
Now we may define the derivative
\begin{multline}
\label{dotx}
(x',\varphi)=\int_I f_0(t)\varphi(t)dt+\sum_{k \in \mathbf N_0}(-1)^{k+1}\int_I f_{k+1}(t)\varphi^{(k+1)}(t)dt+\\+\sum_{k \in \mathbf N_0}\sum_{t \in I} c^{k+1}(t)\sigma_t(\varphi^{(k+1)})+z,
\end{multline}
where $z \in \Delta$,
\begin{equation}
\label{dotx2}
z=\sum_{t \in I} c^0(t)\sigma_t(\varphi)+\sum_{t \in I} \left(\int_{t_1}^t f_0(s)ds\right) \sigma_t(\varphi), \quad \varphi \in \mathcal R.
\end{equation}
As follows from Example \ref{impex}, the distribution defined by (\ref{dotx}), (\ref{dotx2}) is indeed in $D(x)$, and $x'=f$. 
From (\ref{xdef}) we have that $x=0$ for $t<t_0$, and the proof is complete.
\end{proof}

Let us consider in the space $\mathcal R^{n\prime}$ the linear differential equation of the general form
\begin{equation}
\label{eq1}
x'-A(t)x=f
\end{equation}
where $A \in \mathbb G^\infty_{n \times n}$, $f \in \mathcal R^{n\prime}$. The linear systems of form (\ref{eq1}) were considered in \cite{Shi} in the classical space of distributions $\mathcal D'$ in assumption that $A \in \mathbb C_{n \times n}^\infty$. 

We call the solution of the differential equation (\ref{eq1}) the distribution $x \in \mathcal R^{n\prime}$ which possesses the derivative $x' \in \mathcal R^{n\prime}$ such that after substitution of $x$ and $x'$ in (\ref{eq1}) equation (\ref{eq1}) becomes the identity
\begin{equation*}
(x'-A(t)x,\varphi)=(f,\varphi), \quad \varphi \in \mathcal R.
\end{equation*}
Let us note that the left-hand side of differential equation (\ref{eq1}) contains the product of a function from $\mathbb G_{n \times n}^\infty$ (generally discontinuous) and a distribution, which is correctly defined in the space $\mathcal R'$, but in general is undefined in the classical space $\mathcal D'$.

Let $t_0 \in I$. In what follows, we suppose that $f=0$ if $t<t_0$. Let us consider Cauchy problem for the differential equation (\ref{eq1}) with the initial condition
\begin{equation}
\label{ic1}
x=0 \text{ if } t<t_0
\end{equation}
(we will consider Cauchy problem with the initial condition of the general form below).

\begin{theorem}
\label{existsteo}
There exists the unique solution $x$ of Cauchy problem {\rm(\ref{eq1}), (\ref{ic1})}. Furthermore,
$x=Xy$,
where $X$ is the fundamental solution of the corresponding homogeneous system, $y$ is the primitive of the distribution $X^{-1}f$ equal to $0$ if $t<t_0$.
\end{theorem}
\begin{proof}
Without loss of generality we conduct our proof for the case $n=1$. If $X$ is the fundamental solution of the corresponding homogeneous system, then, clearly, $X \in \mathbb G^\infty$ and there exists $X^{-1} \in \mathbb G^\infty$, so the product $X^{-1}f \in \mathcal R'$. According to Theorem \ref{primteo} there exists the primitive of $X^{-1}f$ which is equal to zero if $t<t_0$. Let us denote this primitive by $y$. Now let us determine one of the derivatives of $y$ by the equality
\begin{equation*}
(y',\varphi) := (X^{-1}f,\varphi), \quad \varphi \in \mathcal R.
\end{equation*}
We multiply both sides of $y'=X^{-1}f$ by $X$:
\begin{equation*}
(Xy',\varphi)=(f,\varphi), 
\end{equation*}
so
$Xy'=f$. Let us note that according to the definition of the derivative we have $y'=\frac{dy}{dt}+z$, where $z \in \Delta$. Let us show that $Xz \in \Delta$. Indeed, since $X$ is continuous, we have that if
\begin{equation*}
(z,\varphi)=\sum_{t \in I} c_z(t)\sigma_t(\varphi),
\end{equation*}
then
\begin{equation*}
(Xz,\varphi)=(z,X\varphi)=\sum_{t \in I} c_z(t)X(t)\sigma_t(\varphi),
\end{equation*}
i.e., $c_{Xz}(t):= c_z(t)X(t)$, $t \in I$.
Further, let us choose a particular derivative of the distribution $Xy$. We have:
\begin{equation*}
\left(\frac{d}{dt}(Xy),\varphi \right)=-(Xy,\varphi')=-(y,X\varphi')=-(y,(Xy)^{\prime}-X'\varphi)=\left(X\frac{dy}{dt},\varphi\right)+(X'y,\varphi),
\end{equation*}
where $X' \in \mathbb G^\infty$. 
Since $X\frac{dy}{dt}=f-Xz$, and $Xz \in \Delta$, we may define the derivative
\begin{equation*}
(Xy)^{\prime} := \frac{d}{dt}(Xy)+Xz.
\end{equation*}
Then
\begin{equation*}
(Xy)^{\prime}=f+AXy,
\end{equation*}
so $x=Xy$ is the solution of the differential equation (\ref{eq1}). Further, since $y=0$ for $t<t_0$, then $x=0$ for $t<t_0$: if$\varphi \in \mathcal R$, $\supp(\varphi) \subset (a,t_0)$, then $X\varphi \in \mathcal R$, $\supp(X\varphi) \subset (a,t_0)$ and $$(x,\varphi)=(Xy,\varphi)=(y,X\varphi)=0.$$ Consequently, $x$ is the solution of the Cauchy problem (\ref{eq1}), (\ref{ic1}). 

Finally, let us show that $x$ defined above is the only solution of the Cauchy problem.
Indeed, if $x_1$ and $x_2$ are two solutions of (\ref{eq1}), (\ref{ic1}), then $x_1-x_2$ is the solution of Cauchy problem (\ref{e1}), (\ref{ic1}). According to Theorem \ref{teo2} the difference $x_1-x_2$ is the ordinary solution of system (\ref{e1}), which is identically equal to zero according to classical uniqueness theorem due to condition (\ref{ic1}), so $x_1=x_2$.
\end{proof}

Now let us consider Cauchy problem for the differential equation (\ref{eq1}) with the initial condition of the general form
\begin{equation}
\label{ic2}
x=x_0 \text{ if } t<t_0,
\end{equation}
where $x_0 \in \mathbf C^n$. Consider the following differential equation: 
\begin{equation}
\label{eqgen}
y'-A(t)y=f+x_0\delta_{t_0}^\alpha,
\end{equation}
where $\delta_{t_0}^\alpha \in \mathcal R$, $\alpha \in \mathbf C$ is arbitrary. Clearly, we have $f+x_0\delta_{t_0}^\alpha \in \mathcal R^{n\prime}$, so by Theorem \ref{existsteo} the solution $y$ of Cauchy problem (\ref{eqgen}), (\ref{ic1}) exists in $\mathcal R^{n\prime}$ and unique. We call $y$ the solution of Cauchy problem (\ref{eq1}), (\ref{ic2}). 

The next statement shows that this definition agrees with the classical one.

\begin{theorem} \label{teocool} The following statements are true:

1) The solution of the Cauchy problem (\ref{eq1}), (\ref{ic2}) does not depend on choice of $\alpha$ in (\ref{eqgen}). 

2) If $f$ is a regular distribution, then the solution of Cauchy problem (\ref{eq1}), (\ref{ic2}) coincides for $t>t_0$ with the solution of the corresponding ordinary Cauchy problem
\begin{equation*}
x'-A(t)x=f(t), \quad t>t_0, \quad x(t_0+)=x_0.
\end{equation*}

\end{theorem}
\begin{proof}
Let $y$ be the solution of Cauchy problem (\ref{eqgen}), (\ref{ic1}) (i.e., the solution of Cauchy problem (\ref{eq1}), (\ref{ic2}) according to our definition). Evidently, if $y_1$ is the solution of Cauchy problem (\ref{eq1}), (\ref{ic1}), and
$y_2$ is the solution of Cauchy problem for the differential equation
\begin{equation*}
y'_2-A(t)y_2=x_0\delta_{t_0}^\alpha
\end{equation*}
with initial condition (\ref{ic1}), then by uniqueness of solution $y=y_1+y_2$. Note that $y_1$ does not depend on value of $\alpha \in \mathbf C$, so in order to complete the proof it suffices to show that $y_2$ does not depend on $\alpha$. By Theorem \ref{existsteo} $y_2=Xy_3$, where $X$ is the fundamental solution of the corresponding homogeneous system, and $y_3$ is the primitive of distribution $X^{-1}x_0\delta_{t_0}^\alpha$. Clearly, $X^{-1}$ is continuous on $I$, hence (without loss of generality, $n=1$) $$(X^{-1}x_0\delta_{t_0}^\alpha,\varphi)=(\delta_{t_0}^\alpha,X^{-1}x_0\varphi)=X^{-1}(t_0)x_0\varphi(t_0),$$ and the proof of the first statement is complete.

In order to prove the second statement, first observe that given regular $f$, the solution $x$ of Cauchy problem (\ref{eq1}), (\ref{ic2}) is a regular distribution. Second, since $x_0\delta_{t_0}^\alpha$ coincides with the zero distribution on $I \setminus \{t_0\}$, $x$ coincides with solutions $x_1$, $x_2$ of differential equation
$x'=A(t)x+f(t)$
on $(a,t_0)$ and $(t_0,b)$, respectively. By definition, $x=x_1=0$ on $(a,t_0)$. Further, since $A$ is in $\mathbb G^\infty_{n \times n}$ and, thus, locally-summable on $I$ together with $f$, there exists the limit $x(t_0+)$. Finally, since for every $\alpha$ delta-function $\delta_{t_0}^\alpha$ is the derivative in $\mathcal R'$ of Heaviside function $\theta_{t_0}$ discontinuous at $t_0$, we have that $$x(t_0+)-x(t_0-)=x_0\bigl(\theta_{t_0}(t_0+)-\theta_{t_0}(t_0-)\bigr),$$ so $x(t_0+)=x_0$.
\end{proof}

\begin{theorem}
\label{condep}
The solution of the Cauchy problem (\ref{eq1}), (\ref{ic2}) depends continuously on $f \in \mathcal R^{n\prime}$ and $x_0 \in \mathbf C^n$.
\end{theorem}
\begin{proof}
By definition, the solution of Cauchy problem (\ref{eq1}), (\ref{ic2}) is the solution of Cauchy problem (\ref{eqgen}), (\ref{ic1}). Note that the continuous dependence of the solution on the initial value $x_0$ will follow from the continuous dependence on $f$, so we may restrict ourselves to the proof of the first statement. Now, as follows from the decomposition obtained in the proofs of Theorem \ref{primteo}  (see (\ref{xdef}), (\ref{dotx})), given distribution in $\mathcal R'$, its primitive, equal to $0$ for $t<t_0$, depends continuously on this distribution. Thus, since all operations in $\mathcal R'$ which arise in construction of the solution of the Cauchy problem (\ref{eq1}), (\ref{ic1}) (in particular, (\ref{eqgen}), (\ref{ic1})) in Theorem \ref{existsteo} are continuous (including the operation of multiplication, which is continuous by Theorem \ref{thm3}), the solution depends continuously on $f \in \mathcal R'$.
\end{proof}

\subsection{Linear differential equations of higher orders}
Let us consider in $\mathcal R^{n \times n \prime}$ the following linear differential equation of order $m$:
\begin{equation}
\label{meq}
X^{(m)}-A_{m-1}X^{(m-1)}-\dots-A_0 X=F,
\end{equation}
where $A_i \in \mathbb G^\infty_{n \times n}$, ($0 \leqslant i \leqslant m$), $F \in \mathcal R^{n \times n\prime}$. 

The solution of differential equation (\ref{meq}) is the distribution $X \in \mathcal R^{n \times n\prime}$ which possesses its $m$-th derivative $X^{(m)}$ such that after substitution of $X$, $X'$, \dots, $X^{(m)}$ in (\ref{meq}) we obtain the equality
\begin{equation*}
(X^{(m)}-A_{m-1}X^{(m-1)}-\dots-A_0 X,\varphi)=(F,\varphi)
\end{equation*}
for all $\varphi \in \mathcal R$. 
This definition is correct according to Theorem \ref{importantremark}: it is sufficient to assume that only $X$ and $X^{(m)}$ are specified in order to have all intermediate derivatives $X^{(k)}$, $1 \leqslant k \leqslant m-1$, uniquely determined.

First, let us consider the Cauchy problem for equation (\ref{meq}) with the homogeneous initial conditions
\begin{equation}
\label{mic4}
X^{(k)}=0, \quad t<t_0, \quad 0 \leqslant k \leqslant m-1,
\end{equation}
where $t_0 \in I$. \\ [3mm]
\textbf{Theorem 3.3'.~}There exists the unique solution $X$ of Cauchy problem (\ref{meq}), (\ref{mic4}). Furthermore,
$X=\sum_{i=1}^m X_iT_i ,$
where $\{X_i\}_{i=1}^m$ is the fundamental system of solutions of the corresponding homogeneous differential equation, $T_i$ is the primitive of distribution $Z_iF$ equal to $0$ if $t<t_0$, and $Z_i$ is the $i$-th element of the first column of $R^{-1}$, where
$$R=\left(\begin{array}{ccc}X_1&\dots&X_m\\ & \dots & \\ X_1^{(m-1)}& \dots & X_m^{(m-1)}  \end{array}  \right).$$
\begin{proof}
Following standard scheme, we can make a substitution $Y_k=X^{(k-1)}$, $1 \leqslant k \leqslant m$, to reduce differential equation (\ref{meq}) to the linear system of the form (\ref{eq1}),
\begin{equation}
\label{meq2}
\left(\begin{array}{l}Y_m'\\Y_{m-1}'\\ \dots \\ Y_1' \end{array} \right)
=\left(\begin{array}{rrrrr} A_{m-1}&A_{m-2}&\dots&A_1&-A_0 \\ 1&0&\dots&0&0 \\ && \dots && \\ 0&0&\dots&1&0 \end{array} \right)\left(\begin{array}{l}Y_m\\Y_{m-1}\\ \dots \\ Y_1 \end{array} \right)+\left(\begin{array}{l}F\\0\\ \dots \\ 0 \end{array} \right)
\end{equation}
As it immediately follows from the definition, every solution of differential equation (\ref{meq}) is the solution of system (\ref{meq2}) (see above), and vice versa. The homogeneous initial conditions (\ref{mic4}), being rewritten for system (\ref{meq2}), coincide with the homogeneous initial conditions 
\begin{equation}
\label{moreic}
Y_k=0, \quad t<t_0, \quad 1 \leqslant k \leqslant m, 
\end{equation}
so it suffices to apply Theorem \ref{existsteo} to prove the existence and uniqueness of the solution of Cauchy problem {\rm(\ref{meq}), (\ref{mic4})}. Further, if $\{X_i\}_{i=1}^m$ is the fundamental system of solutions of the corresponding homogeneous linear differential equation, then $R$ is the fundamental matrix of the homogeneous system corresponding to (\ref{meq2}). According to Theorem \ref{existsteo} the solution of Cauchy problem (\ref{meq2}), (\ref{moreic}) admits representation in the form of the product $RT$, where $T=(T_i)_{i=1}^m$ is the primitive of the distribution $R^{-1}(F,0,\dots,0)^{\top}$ which is equal to $0$ for $t<t_0$. Since $X=Y_1$, we have that $X=\sum_{i=1}^m X_i T_i$. The representation of $T_i$ is straightforward. The proof is complete.
\end{proof}

Now consider Cauchy problem for the differential equation (\ref{meq}) with the initial conditions of the general form, that is,
\begin{equation}
\label{mic}
X^{(k)}=X_k, \quad t<t_0, \quad 0 \leqslant k \leqslant m-1,
\end{equation}
where $X_k \in \mathbf C^{n \times n}$. We need to define the solution of Cauchy problem (\ref{meq}), (\ref{mic}).
For this purpose we consider the Cauchy problem for the differential equation of the form (\ref{meq}),
\begin{equation}
\label{meq3}
\sum_{k=0}^m A_k \left(Z^{(k)}-\sum_{j=0}^{k-1} X_j \delta_{t_0}^{(k-j-1)\alpha} \right)=F,
\end{equation}
where $\alpha_{} \in \mathbf C$ is arbitrary, $A_m$ is defined to be the identity matrix, with the homogeneous initial conditions
\begin{equation}
\label{mic2}
Z^{(k)}=0, \quad t<t_0, \quad 0 \leqslant k \leqslant m-1,
\end{equation}
that was already considered above.
So, if $Z$ is the solution of problem (\ref{meq3}), (\ref{mic2}), which always exists and unique in virtue of Theorem 3.3', then we call $Z$ the solution of the Cauchy problem (\ref{meq}), (\ref{mic}). 
The following result shows that this definition is quite natural. \\ [3mm]
\textbf{Theorem 3.4'. }
The following statements are true:

1) The solution of the Cauchy problem (\ref{meq}), (\ref{mic}) does not depend on $\alpha$ in (\ref{meq3}). 

2) If $F$ is a regular distribution, then the solution of the Cauchy problem (\ref{meq}), (\ref{mic}) coincides for $t>t_0$ with the solution of the corresponding ordinary Cauchy problem
\begin{equation*}
X^{(m)}-A_{m-1}X^{(m-1)}-\dots-A_0 X=F, \quad X^{(k)}(t_0+)=X_k, \quad 0 \leqslant k \leqslant m-1.
\end{equation*}
\begin{proof}
Let us consider Cauchy problem for the linear system
\begin{equation}
\label{mmm}
\left(\begin{array}{l}Y_m'\\Y_{m-1}'\\ \dots \\ Y_1' \end{array} \right)
=\left(\begin{array}{rrrr} A_{m-1}&\dots&A_1&A_0 \\ 1&\dots&0&0 \\ & \dots && \\ 0&\dots&1&0 \end{array} \right)\left(\begin{array}{l}Y_m\\Y_{m-1}\\ \dots \\ Y_1 \end{array} \right)+\left(\begin{array}{l}F+X_{m-1}\delta_{t_0}^\alpha\\X_{m-2}\delta_{t_0}^\alpha\\ \dots \\ X_0\delta_{t_0}^\alpha \end{array} \right),
\end{equation}
with the initial conditions as in (\ref{moreic}). As it follows from Theorems \ref{existsteo} and 3.3', the solutions of Cauchy problems (\ref{mmm}), (\ref{moreic}) and (\ref{meq}), (\ref{mic}), respectively, exist and unique. So, in order to prove this theorem we have to show that $Y_1=Z$, the rest will follow from Theorem \ref{teocool} being applied to system (\ref{mmm}). Suppose that $Z$ is the solution of (\ref{meq}), (\ref{mic}). Let us put $$Y_k=Z^{(k-1)}-\sum_{j=0}^{k-2} X_j \delta_{t_0}^{(k-j-2)\alpha}, \quad 2 \leqslant k \leqslant m, \quad Y_1=Z$$
(note that $Z^{(k)}$, $1 \leqslant k \leqslant m$, are uniquely determined, see above).
Let us show that $(Y_k)_{k=1}^m$ is the solution of the Cauchy problem (\ref{mmm}), (\ref{moreic}). Clearly, the initial conditions are satisfied. Let us specify the derivatives $Y_k'$. We define
\begin{equation*}
Y_k':=Z^{(k)}-\sum_{j=0}^{k-2} X_j \delta_{t_0}^{(k-j-1)\alpha}, \quad 1 \leqslant k \leqslant m.
\end{equation*}
Observe that $Y_k'$ is indeed the derivative of $Y_k$. Now for every $1 \leqslant k \leqslant m-1$ we have
\begin{equation*}
Y_k'-Y_{k+1}=X_{k-1}\delta_{t_0}^\alpha.
\end{equation*}
Further,
\begin{multline}
\notag
Y_m'-\sum_{k=0}^{m-1} A_{k} Y_{k+1}=Z^{(m)}-\sum_{j=0}^{m-2} X_j \delta_{t_0}^{(m-j-1)\alpha}-\sum_{k=0}^{m-1} A_{k} \left(Z^{(k)}-\sum_{j=0}^{k-1} X_j \delta_{t_0}^{(k-j-1)\alpha} \right)=\\
=\sum_{k=0}^m A_k \left(Z^{(k)}-\sum_{j=0}^{k-1} X_j \delta_{t_0}^{(k-j-1)\alpha} \right)+X_{m-1}\delta_{t_0}^\alpha=F+X_{m-1}\delta_{t_0}^\alpha,
\end{multline}
so $(Y_k)_{k=1}^m$ is the solution of the Cauchy problem (\ref{mmm}), (\ref{moreic}), as required.
\end{proof}

The next statement is an analogue of Theorem 3.5. \\ [2mm]
\textbf{Theorem 3.5'. }The solution of the Cauchy problem (\ref{meq}), (\ref{mic}) depends continuously on $F \in \mathcal R^{n\times n\prime}$ and $X_k \in \mathbf C^{n \times n}$.

\begin{proof}
The proof follows from the possibility of reduction of linear differential equation (\ref{meq}) to linear system (\ref{mmm}) as in the proof of Theorem 3.4' and from the statement of Theorem \ref{condep} applied to system (\ref{mmm}).
\end{proof}

\begin{proof}[Proof of Example \ref{example1}] Let us find the solution
of Cauchy problem (\ref{pr2}) (the solutions of Cauchy problems (\ref{pr1}) and (\ref{pr3}) can be found similarly). According to Theorem \ref{existsteo}, if $x$ is the solution of Cauchy problem (\ref{pr1}), then $x$ admits the representation $$x=e^{a\zeta_\tau}y,$$ where $e^{a\zeta_\tau(t)}$ is the fundamental solution of the corresponding homogeneous equation, and $y$ is the primitive of the distribution $e^{-a\zeta_\tau}b\delta_\tau^{\prime\alpha}$ which is equal to $0$ for $t<t_0$, so
\begin{multline}
\notag
(e^{-a\zeta_\tau}b\delta_\tau^{\prime\alpha},\varphi)=(\delta_\tau^{\prime\alpha},be^{-a\zeta_\tau}\varphi)=\\=-\alpha b \bigl(-ae^{-\zeta_\tau(\tau+)}\theta_\tau(\tau+)\varphi(\tau+)+e^{-\zeta_\tau(\tau+)}\varphi'(\tau+)\bigr)-\\-(1-\alpha) b \bigl(-ae^{-\zeta_\tau(\tau-)}\theta_\tau(\tau-)\varphi(\tau-)+e^{-\zeta_\tau(\tau-)}\varphi'(\tau-)\bigr)=\\=\alpha b (a\varphi(\tau+)-\varphi'(\tau+))-(1-\alpha)b\varphi'(\tau-),
\end{multline}
so $e^{-\zeta_\tau}b\delta_\tau^{\prime\alpha}=ab\alpha  \delta_\tau^+-b\delta_\tau^{\prime\alpha}$, and the required primitive is $ab\alpha \theta_\tau-b\delta_\tau^\alpha$.
\end{proof}

\begin{proof}[Proof of Example \ref{example2}]
First note that $a \in \mathbb G^\infty$, and $b \in \mathbb L_{\loc}$, so $b \in \mathcal R'$.
As it follows from Theorem \ref{existsteo}, we have to find the primitive of the distribution $e^{-\int_0^t a(a)ds}b''$ which is equal to zero for $t<0$. Due to absolute convergence of the series $\sum_{\gamma \in \Upsilon} b_\gamma$, we have the following chain of equalities:
\begin{multline}
\notag
\left(\exp\left(-\int_0^t a(a)ds\right)b'',\varphi \right)=\\=\left(\exp\left(-\int_0^t a(a)ds\right)\sum_{\gamma \in \Upsilon}b_\gamma\delta_\gamma^{\prime \alpha(\gamma)},\varphi \right)=\sum_{\gamma \in \Upsilon} b_\gamma (\delta_\gamma^{\prime \alpha(\gamma)},e^{-\int_0^t a(a)ds}\varphi)=\\=-\sum_{\gamma \in \Upsilon} b_\gamma e^{-\int_0^\gamma a(a)ds} \bigl(\alpha(\gamma)\varphi'(\gamma+)+\bigl(1-\alpha(\gamma)\bigr)\varphi'(\gamma-)-\\-\alpha(\gamma)a(\gamma+)\varphi(\gamma+)-\bigl(1-\alpha(\gamma)\bigr)a(\gamma-)\varphi(\gamma-) \bigr),
\end{multline}
so the required primitive is
\begin{equation*}
\sum_{\gamma \in \Upsilon}b_\gamma \exp\left(-\int_0^\gamma a(a)ds\right)\left(\bigl(\alpha(\gamma)a(\gamma+)+\bigl(1-\alpha(\gamma)\bigr)a(\gamma-)\bigr)\theta_\gamma-\delta_\gamma^{\alpha(\gamma)} \right),
\end{equation*}
the rest of the proof is straightforward.
\end{proof}

\section{Proofs of Lemmas \ref{lem12}, \ref{lem2} and Theorem \ref{thm3}}

\begin{proof}[Proof of Lemma \ref{lem12}]
By definition, given an absolutely convex set $U \subset \mathcal R$, we have that $U$ is a neighbourhood in $\mathcal R$ if and only if $U \cap \mathbb G^\infty(J)$ is a neighbourhood in $\mathbb G^\infty(J)$, for every $J \in \mathfrak J$. Suppose that $\varphi_k \to \varphi$ in $\mathbb G^\infty$ and there exists $J_0 \in \mathfrak J$ such that $\supp(\varphi_k) \subset J_0$ for all $k \in \mathbf N$. Now, given a neighbourhood of zero $U \subset \mathcal R$ and arbitrary $J \in \mathfrak J$, we have that $$(\varphi_k - \varphi)|_J \in U \cap \mathbb G^\infty(J)$$ starting with certain $k$, since $U \cap \mathbb G^\infty(J)$ is a neighbourhood of zero in $\mathbb G^\infty(J)$, and $(\varphi_k-\varphi)|_J \to 0$ in $\mathbb G^\infty(J)$. The latter follows from the fact that every neighbourhood of zero $U_J \subset \mathbb  G^\infty(J)$ can be extended to a neighbourhood of zero $U_I \subset \mathbb G^\infty$, so that $U_I \cap \mathbb G^\infty(J)=U_J$, while $\varphi_k-\varphi \to 0$ in $\mathbb G^\infty$. By definition, since $J \in \mathfrak J$ was arbitrary, this implies that $\varphi_k \to \varphi$ in $\mathcal R$.

Conversely, suppose that $\varphi_k \to \varphi$ in $\mathcal R$. First, observe that $\varphi_k \to \varphi$ in $\mathbb G^\infty$, since given any neighbourhood $U_I \subset \mathbb G^\infty$, $U_I \cap \mathbb G^\infty(J)$ is a neighbourhood in $ \mathbb G^\infty(J)$ for every $J \in \mathfrak J$, so $U_I$ is a neighbourhood in $\mathcal R$, and, as a result, $\varphi_k-\varphi \in U_I$ starting with certain $k$. Second, suppose that there is no such $J_0 \in \mathfrak J$, and there exists a sequence of subintervals $$\{J_k\}_{k=1}^\infty, \quad J_k \subset J_{k+1}, \quad J_{k+1} \setminus J_k \ne \varnothing, \quad \cap_{k=1}^\infty J_k=I,$$ such that $\varphi_k-\varphi \not \in \mathbb G^\infty (J_k)$ for all $k \in \mathbf N$. The latter may be possible if and only if $\supp(\varphi_k-\varphi) \not \subset \bar{J}_k$. Then we may specify a neighbourhood of zero $U \subset \mathcal R$ such that $\varphi_k-\varphi \not \in U$ for every $k \in \mathbf N$, so $\varphi \not \to \varphi$ in $\mathcal R$, a contradiction. Indeed, let us denote $L_k^-$ and $L_k^+$ the left-hand side and the right-hand side half-open components of $J_{k+1} \setminus J_k$. Assume without loss of generality that $L_k^- \ne \varnothing$. 
Let us define $t_k \in L^-_k$ to be such that $|\varphi_k(t_k-)-\varphi(t_k-)|>0$. By definition, sequence $\{t_k\}_{k=1}^\infty$ tends to the left end-point of the interval $I$. Further, define the monotonically decreasing sequence $\{r_k\}_{k=1}^\infty$ such that $$0<r_k<|\varphi_k(t_k-)-\varphi(t_k-)|, \quad k \in \mathbf N,$$ and $r_k \to 0$. Now let $r \in \mathbb C$ be a function such that $r(t)>0$ for all $t \in I$, and $r(t_k)=r_k$. Finally, we may define the required neighbourhood by $$U=\{\psi\in \mathcal R: |\psi(t)|<r(t), t \in I\}.$$ For any $J \in \mathfrak J$ $\max_{t \in \bar{J}}|r(t)|>0,$ so $U \cap \mathbb G^\infty(J)$ is a neighbourhood in $\mathbb G^\infty(J)$.
\end{proof}

In order to prove Theorem \ref{thm3} we will need the following lemma, whose proof is identical to the proof of the analogous statement for the space $\mathcal D'$ in \cite{Shi}.

\begin{lemma}
\label{rlem}
Given $\{f_k\}_{k=1}^\infty \subset \mathcal R'$, $\{\varphi_k\}_{k=1}^\infty \subset \mathcal R$ such that $\varphi_k \to 0$ in $\mathcal R$ and for every $\varphi \in \mathcal R$ there exists $\lim (f_k,\varphi) \in \mathbf C$, we have that $(f_k,\varphi_k) \to 0.$
\end{lemma}

\begin{proof}
Suppose the contrary. Then we may assume (consider a subsequence, if necessary, ) that there exists $c>0$ such that $|(f_k,\varphi_k)| \geqslant c$, $k \in \mathbf N$. Since $\varphi_k \to 0$ in $\mathcal R$, we may suppose that
$$\|\varphi^{(j)}_k\|_{\mathbb L_{\infty}} \leqslant \frac{1}{4^k}$$ for all $j \leqslant k$. Let us put $\psi_k=2^k\varphi_k$. Then
\begin{equation}
\label{oneineq100}
\|\psi^{(j)}_k\|_{\mathbb L_{\infty}} \leqslant \frac{1}{2^k},
\end{equation}
for all $j \leqslant k$, so $\psi_k \to 0$ in $\mathcal R$, though $$|(f_k,\psi_k)|=2^k|(f_k,\varphi_k)| \geqslant 2^kc \to \infty.$$
Now let us choose $f_{k_1}$, $\psi_{k_1}$ such that
$|(f_{k_1},\psi_{k_1})|>1$. Suppose that $f_{k_j}$, $\psi_{k_j}$ are defined, $1 \leqslant j \leqslant l-1$. Suppose that for every $k \geqslant k'$ we have $$|(f_{k_j},\psi_k)|<\frac{1}{2^{l-j}}$$ for $1 \leqslant j \leqslant l-1$.
Then there exists $k_l \geqslant k'$ such that
\begin{equation}
\label{oneineq101}
|(f_{k_l},\psi_{k_l})|>\sum_{j=1}^{l-1}|(f_{k_l},\zeta_{k_j})|+l
\end{equation}
since $|(f_k,\psi_k)| \to \infty$, we have that $(f_k,\psi_{k_j}) \to 0$ as $k \to \infty$. Suppose that the sequence $\{\psi_{k_l}\}_{l=1}^\infty$ is constructed. Let us define
$\psi=\sum_{j=1}^\infty \psi_{k_j},$
where the series converges due to (\ref{oneineq100}),
so $\psi \in \mathcal R$. Consequently, 
\begin{equation*}
(f_{k_l},\psi)=\sum_{j=1}^{l-1}(f_{k_l},\psi_{k_j})+(f_{k_l},\psi_{k_l})+\sum_{l+1}^\infty (f_{k_l},\psi_{k_j}). 
\end{equation*}
Since (\ref{oneineq101}) and 
\begin{equation*}
\sum_{j=l+1}^\infty (f_{k_l},\psi_{n_j})<\sum_{j=l+1}^\infty \frac{1}{2^{j-l}}=1, 
\end{equation*}
we obtain that $|(f_{k_l},\psi)|>l-1.$ This contradicts to the equality $\lim(f_k,\psi)=(f,\psi)$, where $f=\lim f_k$, so the proof is complete.
\end{proof}
\begin{proof}[Proof of Theorem \ref{thm3}]
Note that $g_k \varphi \to g\varphi$ in $\mathcal R$
for every $\varphi \in \mathcal R(\Omega)$, so 
\begin{multline}
\notag
|(g_kf_k,\varphi)-(gf,\varphi)|=|(f_k,g_k\varphi)-(f,g\varphi)| \leqslant \\
\leqslant  |(f_k,g_k\varphi)-(f_k,g\varphi)|+|(f_k,g\varphi)-(f,g\varphi)|\leqslant\\
\leqslant |(f_k,g_k\varphi-g\varphi)|+|(f_k,g\varphi)-(f,g\varphi)| \to 0
\end{multline}
according to Lemma \ref{rlem} and due to the fact that $f_k \to f$ in $\mathcal R'$.
%where the first summand tends to 0 in virtue of Lemma \ref{dist_complete}, the second summand tends to 0 since $f_k \to f$ in $\mathcal R'(\Omega)$.
\end{proof}

Let $\mathbb {PC} \subset \mathbb G$ be the space of piece-wise constant functions $I \mapsto \mathbf C$. 

\begin{lemma}[\cite{Dieu}]
\label{dieulem}
$\mathbb {PC}$ is dense in $\mathbb G$.
\end{lemma}

\begin{lemma}
\label{lem4} 
The subspace $\mathbb F$ is dense in $\mathbb G^\infty$. Furthermore, for every
$f \in \mathbb G^\infty$ there exists a subsequence $\{f_l\}_{l=1}^\infty \subset \mathbb F$ such that $f_l \to f$
in $\mathbb G^\infty$, $f_l^{(l+1)} \in \mathbb C^{\infty}$, $l \in \mathbf N$, and $$\|f_{l}^{(j)}-f^{(j)}\|_{\mathbb L_{\infty}}<l^{-1}$$ for all $l \geqslant j$, $j \in \mathbf N_0$.
\end{lemma}

\begin{proof}
Let $f \in \mathbb G$. 
First, for every $k \in \mathbf N_0$ we define  the subsequence $\{f_{lk}\}_{l=1}^\infty \subset \mathbb F$ such that $$f_{lk}^{(j)} \to f^{(j)}$$ in $\mathbb L_{\infty}$ for all $0 \leqslant j \leqslant k$. Let us note that the algebra $\mathbb{PC} \subset \mathbb F$ is closed under the operations of integration and differentiation almost everywhere. According to Lemma \ref{dieulem} there exists a subsequence $\{p_{l}^k\}_{l=1}^\infty \subset \mathbb{PC}$ such that  $p_l^k \to f^{(k)}$ in $\mathbb L_{\infty}$, so
$$\int_a^t p_l^k(s)ds \to \int_a^t f^{(k)}(s)ds$$
in $\mathbb L_{\infty}$. Let
\begin{equation*}
p_l^{k-1}(t)=\int_a^t p_l^k(s)ds+q_l^{k-1}(t), 
\quad \{q_l^{k-1}\}_{l=1}^\infty \subset \mathbb P, \quad q_l^{k-1} \to f^{k-1}-\int_a^t f^{(k)}(s)ds
\end{equation*}
in $\mathbb L_{\infty}$, so $p_l^{k-1} \to f^{k-1}$ in $\mathbb L_{\infty}$, $(p_l^{k-1})^{\prime}=p_l^k.$ Further, let us define $\{p_l^{k-2}\}_{l=1}^\infty$,
\begin{equation*}
p_l^{k-2}(t)=\int_a^t p_l^{k-1}(s)ds+q_l^{k-2}(t), 
\quad \{q_l^{k-2}\}_{l=1}^\infty \subset \mathbb P, \quad q_l^{k-2} \to f^{k-2}-\int_a^t f^{(k-1)}(s)ds
\end{equation*}
in $\mathbb L_{\infty}$,
so $p_l^{k-2} \to f^{k-2}$ in $\mathbb L_{\infty}$, $(p_l^{k-2})^{\prime}=p_l^{k-1}.$ 
We may continue this process, and find $p_l^0 \in \mathbb F$. 

Let us define $f_{lk} := p_l^0$. Then $f^{(j)}_{lk}=p_l^j \in \mathbb F$, and $$f^{(j)}_{lk} \to f^{(j)}$$ in $\mathbb L_{\infty}$, where $0 \leqslant j \leqslant k$. We may assume that
\begin{equation*}
\|f_{lk}^{(j)}-f^{(j)}\|_{\mathbb L_{\infty}}<\frac{1}{l}, \quad l \in \mathbf N, \quad 0 \leqslant j \leqslant k.
\end{equation*}
Now define $f_l:=f_{ll} \in \mathbb F$ ($l \in \mathbf N$). Then
$\|f_{l}^{(j)}-f^{(j)}\|_{\mathbb L_{\infty}}<\frac{1}{l},$ $l \geqslant j$, $j \in \mathbf N_0$,
so $f_l \to f$ in $\mathbb G^\infty$.

The sequence $\{f_{l}\}_{l=1}^\infty$ constructed above is the one required, since $f^{(l+1)}_{l} \equiv 0.$
\end{proof}

\begin{proof}[Proof of Lemma \ref{lem2}] Let $\varphi \in \mathcal R$. According to Lemma \ref{lem4} there exists $\{f_k\}_{k=1}^\infty \subset \mathbb{PC}$ such that $f_k \to \varphi$ in $\mathbb G^\infty$. 
Clearly, there exists a test function $\xi \in \mathcal D$ such that $\xi \equiv 1$ in certain open neighbourhood of $\supp(\varphi)$, $\xi \equiv 0$ in certain (larger) open neighbourhood of $\supp(\varphi)$, which possesses the compact closure in $I$. Then $\xi f_k \in \mathcal R_F$, $\xi\varphi=\varphi$, and $\xi f_k \to \xi \varphi$ in $\mathcal R$. Since $\varphi \in \mathcal R$ was chosen arbitrarily, the subspace $\mathcal R_F$ is dense in $\mathcal R$.

If $\{f_k\}_{k=1}^\infty$ is chosen as in the statement of Lemma \ref{lem4}, and we have chosen $\xi \equiv 1$ in certain neighbourhood $t$, then $\{\xi f_k\}_{k=1}^\infty \subset \mathcal R_F$ is the sequence required in the last statement.
\end{proof}

\bibliographystyle{alpha}
\bibliography{kinz}

\begin{thebibliography}{{Gam}69}

\bibitem[BK07]{BK}
A.~{Brudnyi} and D.~{Kinzebulatov}.
\newblock On uniform subalgebras of ${L}^\infty$ on the unit circle generated
  by almost periodic functions.
\newblock {\em Algebra and Analysis}, 19:1--33, 2007.

\bibitem[{Der}86]{Der1}
V.~{Derr}.
\newblock Quasidifferential equations: on the definition of solution of the
  ordinary equation with distributions in coefficients.
\newblock {\em Deponed in VINITI}, 4867-86, 1986.

\bibitem[{Der}02]{Der2}
V.~{Derr}.
\newblock A generalization of \uppercase{R}iemann-\uppercase{S}tieltjes
  integral.
\newblock {\em Func. Diff. Equations}, 1:325--341, 2002.

\bibitem[{Die}69]{Dieu}
J.~{Dieudonne}.
\newblock {\em Foundations of Modern Analysis}.
\newblock Academic Press, 1969.

\bibitem[DK06]{DerKin3}
V.~{Derr} and D.~{Kinzebulatov}.
\newblock Distributions with dynamic test functions and multiplication by
  discontinuous functions.
\newblock {\em Preprint, arXiv:math.CA/0603351}, 2006.

\bibitem[DK07a]{DK2}
V.~{Derr} and D.~{Kinzebulatov}.
\newblock Dynamical generalized functions and the multiplication problem.
\newblock {\em Russian Mathematics}, 51:32--43, 2007.

\bibitem[DK07b]{DK}
V.~{Derr} and D.~{Kinzebulatov}.
\newblock The space of distributions with discontinuous test functions and a
  family of zero-sum games with discontinuous payoffs.
\newblock {\em Preprint, arXiv:math.FA/0606126}, 2007.

\bibitem[{Fil}88]{Fil}
A.F. {Filippov}.
\newblock {\em Differential equations with discontinuous right-hand sides}.
\newblock Kluwer Acad. Publ., 1988.

\bibitem[{Gam}69]{Gam}
T.W. {Gamelin}.
\newblock {\em Uniform algebras}.
\newblock Prentice Hall, 1969.

\bibitem[GS64]{Gelf}
I.M. {Gelfand} and G.E. {Shilov}.
\newblock {\em Generalized functions. Volume 1. Properties and operations}.
\newblock Academic Press, 1964.

\bibitem[KA82]{Kan}
L.V. {Kantorovich} and G.P. {Akilov}.
\newblock {\em Functional Analysis}.
\newblock Pergamon Press, 1982.

\bibitem[KB98]{Kur2}
P.~{Kurasov} and J.~{Boman}.
\newblock Finite rank singular perturbations and distributions with
  discontinuous test functions.
\newblock {\em Proc. Amer. Math. Soc.}, 126:1673--1683, 1998.

\bibitem[{Kin}07]{Kinz}
D.~{Kinzebulatov}.
\newblock Systems with distributions and viability theorem.
\newblock {\em J. Math. Anal. Appl.}, 331:1046--1067, 2007.

\bibitem[{Kur}58]{Kurz2}
J.~{Kurzweil}.
\newblock Generalized ordinary differential equations.
\newblock {\em Chezhosl. Math. Journal}, 8, 1958.

\bibitem[{Kur}59]{Kurz1}
J.~{Kurzweil}.
\newblock Linear differential equations with distributions as coefficients.
\newblock {\em Bull. Acad. Polon. Sci. Ser. Math}, 9, 1959.

\bibitem[{Kur}96]{Kur}
P.~{Kurasov}.
\newblock Distribution theory with discontinuous test functions and
  differential operators with generalized coefficients.
\newblock {\em J. Math. Anal. Appl.}, 201:297--323, 1996.

\bibitem[{Mil}93]{Mil}
B.M. {Miller}.
\newblock The method of discontinuous time substitution in problems of the
  control for impulse and discrete-continuous systems.
\newblock {\em Automat. Remote Control}, 54:1727--1750, 1993.

\bibitem[PT93]{Tvr2}
M.~{Pelant} and M.~{Tvrdy}.
\newblock Linear distributional differential equations in the space of
  regulated functions.
\newblock {\em Mathematica Bohemica}, 118:379--400, 1993.

\bibitem[RR64]{Rob}
A.~{Robertson} and W.~{Robertson}.
\newblock {\em Topological vector spaces}.
\newblock Cambridge Univ. Press, 1964.

\bibitem[{Sar}94]{Sar}
C.O.R. {Sarrico}.
\newblock Some distributional products with relativistic invariance.
\newblock {\em Portugalie Math}, 51:283--290, 1994.

\bibitem[{Sar}95]{Sar2}
C.O.R. {Sarrico}.
\newblock The linear {C}auchy problem for a class of differential equations
  with distributional coefficients.
\newblock {\em Portugalie Math}, 52:379--390, 1995.

\bibitem[{Sar}03]{Sar3}
C.O.R. {Sarrico}.
\newblock Distributional products and global solutions for nonconservative
  inviscid {B}urgers equation.
\newblock {\em J. Math. Anal. Appl.}, 281:641--656, 2003.

\bibitem[{Sch}50]{Scw}
L.~{Schwartz}.
\newblock {\em Th\'eorie des distributions. {Tome} 1}.
\newblock Paris, 1950.

\bibitem[{Shi}84]{Shi}
G.E. {Shilov}.
\newblock {\em Mathematical Analysis. \uppercase{T}he Second Special Course}.
\newblock Moscow Univ. Press, 1984.

\bibitem[SZ97]{Ses}
A.N. {Sesekin} and S.T. {Zavalishin}.
\newblock {\em Dynamic {I}mpulse {S}ystems: {T}heory and {A}pplications.}
\newblock Dordrecht etc.: Kluwer Acad. Publ., 1997.

\bibitem[{Tvr}02]{Tvr1}
M.~{Tvrdy}.
\newblock Differential and integral equations with regulated solutions.
\newblock {\em Mem. Differential Equations Math. Phys.}, 25:1--104, 2002.

\end{thebibliography}

\end{document}